\documentclass[12pt,leqno]{article}

\usepackage{latexsym}
\usepackage{amsmath,amsfonts, theorem}
\usepackage[all]{xy}

\def\B{{\mathcal B}}

\def\C{{\mathbb C}}
\def\F{{\mathbb F}}
\def\N{{\mathbb N}}
\def\Z{{\mathbb Z}}
\def\R{{\mathbb R}}
\def\Q{{\mathbb Q}}
\def\T{{\mathbb T}}
\def\rank{\mathop{\rm rank}}
\def\card{\mathop{\rm card}}
\def\Ad{\mathop{\rm Ad}\,}
\def\id{\mathop{\rm id}}
\def\ita{{\rm it}}

\def\Tr{{\rm Tr}}
\def\Aut{\mathop{\rm Aut}}

\def\sa{{\rm sa}}

\newtheorem{theorem}{Theorem.}[section]
\newtheorem{lemma}[theorem]{Lemma}
\newtheorem{corollary}[theorem]{Corollary}
\newtheorem{definition}[theorem]{Definition}
\newtheorem{proposition}[theorem]{Proposition}

\newtheorem{notation}[theorem]{Notation}
\newtheorem{remark}[theorem]{Remark}

\newenvironment{proof}{\vskip.3cm\noindent{\it Proof.}$\;$}

\numberwithin{equation}{section}

\title{A survey of noncommutative dynamical entropy}
\author{
Erling St\o{}rmer\\
{\small Department of Mathematics, University of Oslo,}\\
{\small P.O. Box 1053, Blindern, 0316 Oslo, Norway}}
\date{}

\begin{document}

\maketitle

\section{Introduction}
With the success of entropy in classical ergodic theory it became a
natural problem to extend the entropy concept to operator algebras. In
the classical case we are given a probability space $(X,B,\mu)$
together with a measure preserving nonsingular transformation $T$ of
$X$. If
$P=\{P_1,\ldots,P_k\}$ is a partition of $X$ by sets in
$\B$ then the entropy of $P$ is
\begin{equation}\label{e1.1}
H(P)=\sum_{i=1}^k \eta(\mu(P_i))\;,
\end{equation}
where $\eta$ is the real function on $[0,\infty)$ defined by
$\eta(0)=0$, $\eta(t)=-t\log t$. One shows that the limit
\begin{equation}\label{e1.2}
H(P,T)=\lim_{n\to\infty}\frac{1}{n}H\Big( \bigvee_0^{n-1}T^{-i}P\Big)
\end{equation}
exists and define the entropy of $T$ by
\begin{equation}\label{e1.3}
H(T)=\sup_P H(P,T)\;,
\end{equation}
where the sup is taken over all finite partitions $P$ as above. Since
$T$ defines an automorphism $\alpha_T$ of $L^\infty(X,\B,\mu)$ by
$\alpha_T(f)(x)=f(T^{-1}(x))$, $x\in X$, the entropy definition
immediately extends to an entropy $H(\alpha_T)$, where in the
definition we replace $P$ by the finite dimensional subalgebra of
$L^\infty(X,\B,\mu)$ spanned by the characteristic functions
$\chi_{P_i}$.

If we want to extend this definition to the noncommutative setting the
obvious first try is first to define the entropy $H_\varphi(N)$ of a
finite dimensional algebra with respect to a state $\varphi$. Then by
analogy with (\ref{e1.2}) and (\ref{e1.3}) if $\alpha$ is a
$\varphi$-invariant state, to consider $\frac{1}{n}H_\varphi\big(
\bigvee\limits_0^{n-1}\alpha^i(N)\big)$, where $\bigvee\limits_{i=1}^r
A_i$ stands for the von Neumann algebra generated by the algebras
$A_1,\ldots,A_r$.
This approach does not work, because the C*-algebra
generated by two finite dimensional C*-algebras need not be finite
dimensional. There have been several approaches to circumvent this
difficulty. We shall consider the two we consider most successful.
The
first was initiated by Connes and St\o{}rmer \cite{C-S} and consisted
of defining a function $H(N_1,\ldots,N_k)$ on finite families
$N_1,\ldots,N_k$ of finite dimensional subalgebras of a von Neumann
algebra with a normal tracial state, which satisfies many of the same
properties as the entropy function $H(N_1\vee\cdots\vee N_k)$ in the
abelian case.
This
was possible due to the beautiful properties of relative entropy of
states and the function $\eta(t)$. Later on Connes, Narnhofer and
Thirring \cite{CNT} extended this definition to entropy with respect
to invariant states on C*-algebras.

The second approach due to Voiculescu \cite{V} is a refinement of the
mean entropy described above. Instead of starting with a finite
dimensional subalgebra $N$ and looking at
$\bigvee\limits_0^{n-1}\alpha^i(N)$, he considered finite subsets
$\omega$ of the von Neumann algebra and then looked for finite
dimensional subalgebras which approximately contained
$\bigcup\limits_0^{n-1}\alpha^i(\omega)$.
Then the
entropy, or the rank, of this algebra was used in the definition.
There are several variations of this definition depending on how the
approximation is taken. They all majorize the C-S or CNT-entropies
indicated in the previous paragraph.

The two definitions behave quite differently with respect to tensor
products. The CNT-entropy is superadditive, i.e.
$h(\alpha\otimes \beta)\geq h(\alpha)+h(\beta)$, while those of
Voiculescu are subadditive. In many cases the two entropies coincide,
so that the tensor product formula
$h(\alpha \otimes \beta)=h(\alpha)+h(\beta)$ holds.

Having the different definitions of entropy the natural question is:
what do they tell us about the automorphism?
In the
classical cases in addition to being a good conjugacy invariant,
entropy roughly measures the amount of ergodicity of the transformation
and how fast and how far finite dimensional subalgebras are moved with
increasing powers of the transformation.
In the
noncommutative situation much the same is true, except for one major
difference. The entropy also measures the amount of commutativity
between finite dimensional subalgebras and their images under the
action. Thus in highly noncommutative cases like infinite free products
of an algebra with itself and the shift, the entropy is zero even
though the shift is extremely ergodic. On the other hand, shifts on
infinite tensor products behave like classical shifts.

The aim of these notes is to describe all the above in more detail
together with the most studied examples. They will usually be
introduced in places where they illustrate and show applications of the
theory.
We shall
rarely give complete proofs, but will indicate the main ideas in many
cases in order to exhibit the mathematical techniques and ideas
involved. The bibliography is not meant to be complete; we have as a
rule tried to include references to papers directly related to the
text. For other approaches and references see \cite{A-F}, \cite{Hu},
\cite{T}.

The notes are organized as follows

In Section~2 we treat the entropy of Connes and St\o{}rmer on finite von
Neumann algebras. We start with the background on the operator concave
function $\eta(t)$ and relative entropy. Then we define the entropy
function $H(N_1,\ldots,N_k)$ and discuss its properties. After defining
entropy of a trace invariant automorphism and stating its basic
properties we illustrate the results by looking at noncommutative
Bernoulli shifts.

Section~3 is devoted to the extension of Connes, Narnhofer and Thirring
of the results in Section~2 to automorphisms and invariant states of
C*-algebras. They replace the entropy function $H(N_1,\ldots,N_k)$ by a
similar function $H(\gamma_1,\ldots,\gamma_k)$ defined on completely
positive maps $\gamma_1,\ldots,\gamma_k$ from finite dimensional
C*-algebras into the C*-algebra.

In Section~4 we consider the example which has attracted most attention
in the theory, namely quasifree states of the CAR-algebra and invariant
Bogoliubov automorphisms. The formula for the CNT-entropy we shall
discuss, has been gradually extended to more general situations,
starting with \cite{SV} and now being completed in \cite{N}.

Section~5 is devoted to the entropy of Sauvageot and Thouvenot
\cite{S-T}. This entropy is a variation of the CNT-entropy, and they
coincide for nuclear C*-algebras and injective von Neumann
algebras. As an illustration type I algebras are considered in some
detail.

In Section~6 we define and study Voiculescu's approximation entropies
\cite{V} together with Brown's extension \cite{Br1} of topological
entropy to exact C*-algebras. These entropies majorize the CNT-entropy
and comparison of them can yield much information on the C*-dynamical
system under consideration.

Section~7 is devoted to crossed products. If $(A,\varphi,\alpha)$ is a
C*-dynamical system, i.e. $A$ is a C*-algebra, $\varphi$ a state and
$\alpha$ a $\varphi$-invariant automorphism, then $\alpha$ extends to an
inner automorphism $\hat{\alpha}$ of $A\times_\alpha\Z$. Using the
theory from Section~6 we obtain results, even in more general
situations, to the effect that the entropies of $\alpha$ and
$\hat{\alpha}$ are the same. In particular the shift on
$\mathcal{O}_\infty$ has topological entropy zero.

In Section~8 we study the most noncommutative setting, namely shifts on
infinite free products $(\ast A_i,\ast \varphi_i)$ where the $A_i$'s
and the $\varphi_i$'s are equal. These automorphisms are ``extremely''
ergodic, but still their entropies vanish.

Section~9 is devoted to binary shifts on the CAR-algebra, arising from
sequences of 0's and 1's. Different bitstreams give rise to
C*-dynamical systems of quite different nature. The entropies with
respect to the trace are in most computed cases equal to
$\frac{1}{2}\log 2$, but there are examples with entropy zero.

In Section~10 on generators we consider W*-dynamical systems
$(M,\tau,\alpha)$ with $M$ a von Neumann algebra with a faithful normal
tracial state $\tau$, where the entropy is a mean entropy. In such
cases the C-S entropy tends to coincide with one of Voiculescu's
approximation entropies. The concept of generator can be made quite
general, and we get in some cases the analogue of the classical formula
when the entropy of a transformation $T$ is the relative (or
conditional) entropy $H(\bigvee\limits_0^\infty T^{-i}P\mid
\bigvee\limits_1^\infty T^{-i}P)$. Applications are given to
subfactors and the canonical endomorphism $\Gamma$ on the hyperfinite
II$_1$-factor defined by an inclusion of subfactors of finite index.

Finally, Section~11 is devoted to the variational principle. Several of
the well-known results from the classical case and from spin lattice
systems in the C*-algebra formalism of quantum statistical mechanics
are extended to a class of asymptotically abelian C*-algebras.

Acknowledgement. The author is indebted to S.~Neshveyev for several
useful comments.

\section{Entropy in finite von Neumann algebras}

In this section we shall define and sketch the proofs of the main
properties of the entropy function $H(N_1,\ldots,N_k)$ and the
corresponding entropy of a trace invariant automorphism. For this we
need to study the function $\eta(t)=-t\log t$, $t>0$, $\eta(0)=0$, and
relative entropy in some detail. The first goes back to early work on
entropy of states, see \cite{N-U}. Recall that $B(H)$ denotes the
bounded linear operators on a Hilbert space $H$, and $B(H)^+$ the
positive operators in $B(H)$.

\begin{lemma}\label{L2.1}
(i) \ The function $\log t$ is operator increasing on $B(H)^+$, i.e. if
$0\leq x\leq y$ in $B(H)^+$ then $\log x\leq \log y$.

(ii) \ The
function $\eta(t)$ is strictly operator concave on $B(H)^+$, i.e.
\[
\eta\big({\textstyle\frac{1}{2}}(x+y)\big)\geq
{\textstyle\frac{1}{2}}\eta(x)
+{\textstyle\frac{1}{2}}\eta(y)\;,\qquad x,y\in B(H)^+
\]
with equality only if $x=y$.
\end{lemma}

\begin{proof}
For $t\geq 0$
\[
\log t=\int\limits_0^\infty\Big(
\frac{1}{1+x}-\frac{1}{t+x}\Big)dx\;,
\]
providing (i). Multiplying by $t$ we get
\[
\eta(t)=\int\limits_0^\infty\Big( 1-
\frac{t}{1+x}-\frac{x}{t+x}\Big)dx\;.
\]
When we take convex combinations, the first two summands cancel out, so
the lemma follows from the inequality
\[
\big({\textstyle\frac{1}{2}}(z+w)\big)^{-1}\leq
{\textstyle\frac{1}{2}}(z^{-1}+w^{-1})
\]
for positive invertible operators $z$ and $w$.
\hfill
$\Box$
\end{proof}
\bigskip

Instead of using partitions of unity consisting of orthogonal
projections as in the classical case it will be necessary to look at
more general partitions of unity.

\begin{notation} \label{N2.2}
\rm
Let $M$ be a von Neumann algebra, and let $k\in\N$. Then
\begin{eqnarray*}
\lefteqn{S_k=S_k(M)=\{ (x_{i_1,\ldots,i_k}):
    x_{i_1,\ldots,i_k}\in M^+\quad \mbox{and equal to 0}} \\
&&\qquad \mbox{except for a finite number of indices,
$\sum\limits_{i_1\ldots i_k} x_{i_1 \ldots i_k}=1$}\}
\end{eqnarray*}
Let for $j\in\{1,\ldots,k\}$
\[
x_{i_j}^j=
\sum_{i_1,\ldots,i_{j-1},i_{j+1},\ldots,i_k}
x_{i_1\ldots i_k}\;.
\]
One can then show the following inequality \cite{H-S}.

\begin{lemma} \label{L2.3}
Let $M$ be a von Neumann algebra with a normal tracial state $\tau$.
Let $\|x\|_2=\tau(x^\ast x)^{1/2}$ for $x\in M$. Let $(x_{ij})\in S_2$,
$i=1,\ldots,m$, $j=1,\ldots,n$. Then
\[
\sum_i \tau\eta(x_i^1)+
\sum_j \tau\eta(x_j^2)
-\sum_{i,j} \tau\eta(x_{ij})\geq
{\textstyle\frac{1}{2}}\sum_{ij} \|
[ (x_i^1)^{1/2},(x_j^2)^{1/2}] \|_2
\]
where $[a,b]=ab-ba$.
\end{lemma}

In particular the left side of the inequality is nonnegative. This
follows also from the joint convexity of {\em relative entropy},
defined as follows: If $x,y\in M^+$ with $x\leq \lambda y$ for some
$\lambda>0$,
\begin{equation}\label{e2.4}
S(x,y)=\tau(x(\log x-\log y))\;.
\end{equation}

More generally if $M$ is a von Neumann algebra and $\varphi,\psi$
normal states we can define their relative entropy as follows,
\cite{A}, \cite{O-P}. We may assume $\varphi$ and $\psi$ are vector
states $\omega_{\xi_\varphi}$ and $\omega_{\xi_\psi}$ respectively and
for simplicity that $\xi_\varphi$ is separating and cyclic for $M$. We
define
\[
S_{\psi,\varphi}(x\xi_\varphi)=x^\ast \xi_\psi
\]
If $\bar{S}_{\psi\varphi}$ is the closure, the relative modular operator
is
\[
\Delta_{\psi,\varphi}=S_{\psi\varphi}^\ast
\bar{S}_{\psi\varphi}\;.
\]
Then the relative entropy is
\[
S(\varphi,\psi)=-(\log \Delta_{\psi,\varphi}
\xi_\varphi, \xi_\varphi)\;.
\]
There are also integral formulas due to Pusz, Woronowicz and Kosaki
which yield generalizations to C*-algebras, see \cite{O-P}. One can
show that $S$ is jointly convex in $\varphi$ and $\psi$ and
$S(\lambda\varphi,\lambda\psi)=
\lambda S(\varphi,\psi)$. Furthermore, if $x\in M^+$ and $S$ is
given by (\ref{e2.4}) then
\[
S(x,\tau(x))=\tau(x(\log x-\log \tau(x))=
\eta(\tau(x))-\tau \eta(x)\;.
\]
This together with joint convexity of $S$ yields the inequality
\begin{equation}\label{e2.5}
\eta \tau(x+y)-\tau \eta(x+y)\leq
(\eta \tau(x)-\tau \eta(x))
+(\eta \tau(y)-\tau \eta(y))\;.
\end{equation}

If $\varphi$ is a normal state on $M$ then there exists a positive
self-adjoint operator $h_\varphi\in L^1(M,\tau)$ such that
$\varphi(x)=\tau(h_\varphi x)$. Then the relative entropy of $\varphi$
and $\omega$ is given by
\[
S(\varphi,\omega)=S(h_\varphi,h_\omega)=
\varphi(\log h_\varphi-\log h_\omega)
\]
whenever it is defined.

If $N\subset M$ is a von Neumann subalgebra we denote by $E_N$ the
trace invariant conditional expectation of $M$ onto $N$ defined by the
identity
\[
\tau(E_N(x)y)=
\tau(xy)\qquad \mbox{for $\;x\in M$, $y\in N$}\;.
\]
If $\varphi$ and $\omega$ are normal states of $N$ and $M$ respectively
we have, see \cite[Thm.~5.15]{O-P},
\begin{equation}\label{e2.6}
S(\omega,\varphi\circ E_N)=
S(\omega|_N,\varphi)+S(\omega,\omega\circ E_N)\;.
\end{equation}
If $\varphi,\psi,\omega\in M_\ast^+$ and $\omega\leq \psi$ then by
\cite[Cor.~5.12]{O-P}
\begin{equation}\label{e2.7}
S(\varphi,\psi)\leq S(\varphi,\omega)\;.
\end{equation}
After these preliminaries we now define the entropy function
$H(N_1,\ldots,N_k)$.
\end{notation}

\begin{definition} \label{D2.4}
{\rm \cite{C-S}} \
Let $N_1,\ldots,N_k$ be finite dimensional von Neumann subalgebras of
$M$. Then
\[
H(N_1,\ldots,N_k)=
\sup_{x_{i_1\ldots i_k}\in S_k} \Big\{
\sum_{i_1\ldots i_k} \eta\tau(x_{i_1\ldots i_k})
-\sum_{j=1}^k \sum_{i_j} \tau\eta(E_{N_j}
x_{i_j}^j)\Big\}\;.
\]
\end{definition}

The definition can be rewritten in terms of relative entropy as
follows. Write $x_{(i)}$ for $x_{i_1\ldots i_k}$. Let $\tau_{(i)}$ and
$\tau_{i_j}^j$ denote the positive linear functionals
\[
\tau_{(i)}(a)=
\tau(x_{(i)}a),\qquad \tau_{i_j}^j(a)=\tau(x_{i_j}^j a)\;.
\]
Since
\begin{eqnarray*}
S(\tau_{i_j}^j|_{N_j},\tau|_{N_j}) &=&
    \tau (E_{N_j}(x_{i_j}^j)(\log E_{N_j}(x_{i_j}^j)
    -\log E_{N_j}(1))) \\
&=& -\tau\eta(E_{N_j}x_{i_j}^j)\;,
\end{eqnarray*}
the definition of $H$ becomes
\begin{equation}\label{2.8}
H(N_1,\ldots,N_k)
=\sup_{(\tau_{(i)})} \Big\{
\sum_{(i)} \eta\tau_{(i)}(1)
+\sum_{j=1}^k \sum_{i_j} S(\tau_{i_j}^j|_{N_j},\tau|_{N_j})\Big\}\;.
\end{equation}

The main properties of $H$ are summarized in

\begin{theorem} \label{T2.10}
{\rm \cite{C-S}} \
For finite dimensional von Neumann subalgebras $N,N_i,P_j$ of $M$ we
have
\begin{itemize}
\item[{\rm (A)}]
$H(N_1,\ldots,N_k)\leq H(P_1,\ldots,P_k)$ when $N_i\subset P_i$,
$i=1,\ldots,k$.

\item[{\rm (B)}]
$H(N_1,\ldots,N_k,N_{k+1},\ldots,N_p)\leq
       H(N_1,\ldots,N_k)+H(N_{k+1},\ldots,N_p)$

\item[{\rm (C)}]
$N_1,\ldots,N_k\subset N\Rightarrow
       H(N_1,\ldots,N_k,N_{k+1},\ldots,N_p)\leq
       H(N,N_{k+1},\ldots,N_p)$

\item[{\rm (D)}]
For any family of minimal projections of $N,(e_j)_{j\in I}$, such that
$\sum\limits_{j\in I}e_j=1$ we have
$H(N)=\sum\limits_{j\in I}\eta\tau(e_j)$.

\item[{\rm (E)}]
If $P_i$ pairwise commute, $P_i\subset N_i$, and
$\bigvee\limits_{i=1}^k P_i=
\bigvee\limits_{i=1}^k N_i$ then
\[
H(N_1,\ldots,N_k)=H\Big(
\bigvee\limits_{i=1}^k N_i\Big)\;.
\]
\end{itemize}
\end{theorem}

\subsubsection*{Indication of proof}

(A) A variant of Jensen's inequality states that $\eta(E_N(x))\geq
E_N(\eta(x))$ for $x\in M^+$. If $N_i\subset P_i$ then
$\eta(E_{N_i}(x))=\eta(E_{N_i}E_{P_i}(x))\geq
E_{N_i}\eta(E_{P_i}(x))$, hence $\tau\eta(E_{N_i}(x))\geq
\tau\eta(E_{P_i}(x))$, proving (A).
\bigskip

\noindent
(B) \
This is a reduction to subadditivity of $H(P)$ in the classical case.
\bigskip

\noindent
(C) \
This is a consequence of the positivity of the left side of the
inequality in Lemma~\ref{L2.3}.
\bigskip

\noindent
(D) \
The proof of this property is helpful in understanding the need for the
second sum in Definition~\ref{D2.4}. Let $(e_j)_{j\in I}$ be as in
(D). Let
$(x_i)\in S_1$. Since $\tau(x_i)=\tau(E_N x_i)$, we may assume $x_i\in
N$. Thus we have to show
\[
\sum_i \eta\tau(x_i)-\sum_i \tau\eta(x_i)\leq
\sum_{j\in I}\tau\eta(e_j)\;.
\]
By inequality (\ref{e2.5}) we can reduce to the case when each $x_i$ is
of $\;\rank 1$, i.e. $x_i=\lambda_ip_i$ with $\lambda_i>0$, $p_i$ a
minimal projection in $N$. Computing and noting that $\eta(p_i)=0$ we
have
\begin{eqnarray*}
\lefteqn{
   \sum(\eta\tau(x_i)-\tau\eta(x_i))=\sum(\eta\tau(\lambda_i p_i)
     -\tau\eta(\lambda_i p_i))} \\
&&=\sum \lambda_i \eta(\tau(p_i))+\eta(\lambda_i)\tau(p_i)
     -\lambda_i \tau\eta(p_i)-\eta(\lambda_i)\tau(p_i) \\
&&=\sum \lambda_i \eta(\tau(p_i))\;.
\end{eqnarray*}
If we write $N$ as a direct sum of factors, we reduce to the case when
$N\cong M_n(\C)$, $n\in\N$, so that
\[
1=\sum \tau(\lambda_i p_i)=\Big( \sum \lambda_i\Big)\frac{1}{n}\;,
\]
hence $\sum \lambda_i=n$, and so
\[
\sum \lambda_i \eta(\tau(p_i))=
   n\cdot \eta\Big(\frac{1}{n}\Big)=\log n=\sum \eta\tau(e_j)\;.
\]
(E) \
This property shows that the definition of $H$ generalizes the abelian
case. In the proof we can replace each $N_i$ by $P_i$. If $A_i$ is a
masa, i.e. a maximal abelian subalgebra of $P_i$, and
$A=\bigvee\limits_{i=1}^k A_i$ is the masa they generate in
$\bigvee\limits_{i=1}^k P_i$, then by (D), $H(A)=H\Big(
\bigvee\limits_{i=1}^k P_i\Big)$. Thus (E) follows from the abelian
case.
\hfill
$\Box$
\bigskip

The function $(N_1,\ldots,N_k)\to H(N_1,\ldots,N_k)$ is from the above
a function of the sizes of the $N_i$'s together with their relative
positions. It seems to be very difficult to formulate a general theorem
in the converse direction. One simple result follows from (D), namely
if $P\subset N$ and $H(P)=H(N)$ then each masa in $P$ is a masa in $N$,
i.e. $\rank P=\rank N$, where the $\rank$ of $N-\rank N=\dim A$, where
$A$ is a masa in $N$. Note that $\dim N\leq(\rank N)^2$. So far the
only theorem in the literature along the line discussed above is

\begin{theorem} \label{T2.6}
{\rm \cite{H-S}} \
Let $M$ and $\tau$ be as before. Let $N_1,\ldots,N_k$ be finite
dimensional von Neumann subalgebras of $M$, and let
$N=\bigvee\limits_{i=1}^k N_i$. Then the following two conditions are
equivalent.
\medskip

(i)\
$H(N_1,\ldots,N_k)=H(N)$

(ii)\
There exists a masa $A\subset N$ such that
$A=\bigvee\limits_{i=1}^k(A\cap N_i)$.

\medskip

\noindent
In particular, if the above
conditions hold then $N$ is finite dimensional, and
\newline
$\rank N\leq
\prod\limits_{i=1}^k\rank N_i$.
\end{theorem}

Note that the implication (ii)$\Rightarrow$(i) is an easy consequence
of Theorem~2.5. Indeed
\begin{eqnarray*}
H(N) &\geq& H(N_1,\ldots,N_k)\qquad \mbox{by (C)} \\
&\geq& H(A\cap N_1,\ldots,A\cap N_k)\qquad \mbox{by (A)} \\
&=& H\Big( \bigvee_{i=1}^k (A\cap N_i)\Big)\qquad \mbox{by (E)} \\
&=& H(A) \\
&=& H(N)\qquad \mbox{by (D)}
\end{eqnarray*}
For the converse we must attack the definition of $H$,
Definition~\ref{D2.4}, directly. Choose $(x_{(i)})\in S_k$ for which the
right side of Definition~\ref{D2.4} almost takes the value
$H(N_1,\ldots,N_k)$. By using the
$k$-dimensional version of the inequality in Lemma~\ref{L2.3} it follows
that the operators  $x_{i_j}^j$ almost commute for different $j$'s,
and taking limits of such families $(x_{(i)})\in S_k$ we can conclude
that the $x_{i_j}^j$ belong to pairwise commuting algebras $P_j$.
Taking masas $A_j$ in these $P_j$ we get the desired $A$ as
$A=\bigvee\limits_{j=1}^k A_j$.
\hfill
$\Box$
\bigskip

In the classical case two finite dimensional algebras $A$ and $B$
(identified with the partition of unities of their atoms) are said to
be independent if $\mu(fg)=\mu(f)\mu(g)$, $f\in A$, $g\in B$, or
equivalently $H(A\vee B)=H(A)+H(B)$. This equivalence is false in the
noncommutative case. However, we have

\begin{corollary} \label{C2.12}
{\rm \cite{H-S}} \
Let $N_1,\ldots,N_k\subset M$ as before and put
$N=\bigvee\limits_{i=1}^k N_i$. Then the following two conditions are
equivalent.
\begin{itemize}
\item[(i)]
$H(N)=H(N_1,\ldots,N_k)=\sum\limits_{i=1}^k H(N_i)$.

\item[(ii)]
There exists a masa $A\subset N$ such that $A_i=A\cap N_i$ is a masa in
$N_i$ for each $i$, and $A_1,\ldots,A_k$ are independent.
\end{itemize}
\end{corollary}

We shall next consider continuity of $H$.

\begin{definition} \label{D2.8}
If $N,P\subset M$ are finite dimensional subalgebras their {\em
relative entropy\/} is
\[
H(N|P)=\sup_{x\in S_1}\sum_i(\tau\eta(E_P x_i)-\tau\eta(E_N x_i))
\]
\end{definition}

(If we compare with the classical situation we should perhaps rather
have used the name ``conditional entropy''). The following properties
are immediate consequences of Definition~\ref{D2.8}

(F) \
$H(N_1,\ldots,N_k)\leq H(P_1,\ldots,P_k)+\sum\limits_{j=1}^k
    H(N_j|P_j)$.

(G) \
$H(N|Q)\leq H(N|P)+H(P|Q)$.

(H) \
$H(N|P)$ is increasing in $N$ and decreasing in $P$.
\bigskip

If $N\supset P$ the definition makes sense even when $N$ and $P$ are
infinite dimensional, as noted by Pimsner and Popa \cite{P-P}. They
computed $H(N|P)$ in several cases relating it in particular to the
Jones index, see Section~10.

In the classical case the crucial result which makes entropy useful, is
the Kolmogoroff-Sinai theorem, see \cite{Sh} for a natural proof using
continuity of relative entropy. Continuity in our case takes the form
of the following lemma. For $N,P\subset M$ and $\delta>0$ we write
$N\subset^\delta P$ if for each $x\in N$, $\|x\|\leq 1$, there exists
$y\in P$, $\|y\|\leq 1$, such that $\|x-y\|_2<\delta$.

\begin{lemma} \label{L2.9}
{\rm \cite{C-S}} \
Let $M$ and $\tau$ be as before and $n\in\N$, $\varepsilon>0$. Then
there exists $\delta>0$ such that for all pairs of von Neumann
subalgebras $N,P\subset M$ we have:
\[
\dim N=n\;,\qquad N\subset^\delta P\Rightarrow H(N|P)<\varepsilon\;.
\]
\end{lemma}

\begin{definition}\label{D2.15}
{\rm \cite{C-S}} \
Let $\alpha$ be an automorphism of $M$ which is $\tau$-invariant, i.e.
$\tau\circ \alpha=\tau$. If $N$ is a finite dimensional von Neumann
subalgebra of $M$, put
\[
H(N,\alpha)=\lim_{k\to\infty}\frac{1}{k}
H(N,\alpha(N),\ldots,\alpha^{k-1}(N))\;.
\]
\end{definition}

This limit exists by property (B), see \cite[Thm.~4.9]{W}. The {\em
entropy\/} $H_\tau(\alpha)$, or $H(\alpha)$, of $\alpha$ is
\[
H(\alpha)=\sup_N H(N,\alpha)\;,
\]
where the $\sup$ is taken over all $N$ as above. The Kolmogoroff-Sinai
theorem takes the form, see \cite[Thm.~4.22]{W} for the classical
analogue.

\begin{theorem}\label{T2.16}
{\rm \cite{C-S}} \
Let $M$ be hyperfinite, and $\tau$ and $\alpha$ as above. Let
$P_j)_{j\in\N}$ be an increasing sequence of finite dimensional
subalgebras of $M$ with $\Big(\bigcup\limits_{j=1}^\infty
P_j\Big)''=M$. Then
\[
H(\alpha)=\lim_{j\to\infty} H(P_j,\alpha)\;.
\]
\end{theorem}

\begin{proof}
Let $N\subset M$ be finite dimensional and $\varepsilon>0$. By
hypothesis and Lemma~\ref{L2.9} there exists $j\in\N$ such that
$H(N|P_j)<\varepsilon$. Thus by property (F)
\begin{eqnarray*}
H(N,\alpha) &=&
    \lim_k \frac{1}{k}H(N,\alpha(N),\ldots,\alpha^{k-1}(N)) \\
&\leq& \lim_k \frac{1}{k} H(P_j,\alpha(P_j),\ldots,
       \alpha^{k-1}(P_j))+\lim_k\frac{1}{k}
    \sum_{i=0}^{k-1} H(\alpha^i(N),\alpha^i(P_j)) \\
&\leq& H(P_j,\alpha)+\varepsilon\;.
\end{eqnarray*}

\vspace*{-3ex}
\hfill
$\Box$
\end{proof}
\bigskip

It is clear that $H(\alpha)$ is a conjugacy invariant, i.e. if $\gamma$
is an automorphism of $M$ then $H(\gamma\alpha\gamma^{-1})=H(\alpha)$.

In the classical case we have $H(\alpha^p)=|p|H(\alpha)$ for $p\in\Z$.
In our case we have,

\begin{proposition}\label{P2.16}
{\rm \cite{C-S}} \
(i) \
$H(\alpha^p)\leq |p|H(\alpha)$.
\newline
(ii) \
If $M$ is hyperfinite, $H(\alpha^p)=|p|H(\alpha)$.
\end{proposition}

Note also that by property (C) $H(\alpha)$ is monotone, i.e. if
$N\subset R$ is a von Neumann subalgebra such that $\alpha(N)=N$, then
$H(\alpha|_N)\leq H(\alpha)$.

A problem which has attracted much attention in noncommutative entropy
is that of additivity under tensor products. If $(M_i,\alpha_i,\tau_i)$
are W*-dynamic systems like $(M,\tau,\alpha)$ above, $i=1,2$, then the
problem is whether
\[
H_{\tau_1\otimes\tau_2}(\alpha_1\otimes\alpha_2)=
H_{\tau_1}(\alpha_1)+H_{\tau_2}(\alpha_2)\;?
\]
This is well-known in the classical case. In our case we can only
conclude that
\begin{equation}\label{2.18}
H_{\tau_1\otimes\tau_2}(\alpha_1\otimes\alpha_2)\geq
H_{\tau_1}(\alpha_1)+H_{\tau_2}(\alpha_2)\;.
\end{equation}
Indeed, if $N_i\subset M_1$, $P_i\subset M_2$, $i=1,\ldots,k$ are
finite dimensional then there are more families
$(x_{(i)})=(x_{i_1},\ldots,x_{i_k})\in S_k(M_1\otimes M_2)$ then there
are families $(y_{(i)}\otimes z_{(i)})=(y_{i_1\ldots i_k}\otimes
z_{i_1\ldots i_k})$ in $S_k(M_1)\otimes S_k(M_2)$, hence
\[
H_{\tau_1\otimes\tau_2} (N_1\otimes P_1,\ldots,N_k\otimes P_k)
\geq H_{\tau_1}(N_1,\ldots,N_k)+H_{\tau_2}(P_1,\ldots,P_k)\;.
\]

\begin{remark}\label{R2.19}
\rm
{\bf The $n$-shift}\quad
The first nontrivial example that was computed was the entropy of the
$n$-shift. Let $n\in\N$, $M_i=M_n(\C)$, $i\in\Z$, and $\tau_i$ be the
tracial state on $M_i$. Let $B=\bigotimes\limits_{i\in\Z}M_i$,
$\tau=\bigotimes\limits_{i\in\Z}\tau_i$, be the C*-tensor product, and
consider $B$ as a subalgebra of the II$_1$-factor $R$ obtained from the
GNS-representation of $\tau$. Let $\alpha$ be the shift on $B$
identified with its extension to $R$. Let
\[
P_j=\cdot\otimes 1\otimes \bigotimes_{-j}^j M_i\otimes 1\otimes\cdots,
\qquad j\in\N
\]
be the finite tensor product of the $M_i$ from $-j$ to $j$ considered
as a subalgebra of $R$. Let $D_i$ be the diagonal in $M_i$ and
\[
D_{pq}=\cdots\otimes 1\otimes\bigotimes_{-p}^q D_i\otimes\cdots,
\qquad j\in\N\;.
\]
As an illustration of the techniques developed we give two computations
of $H(\alpha)$. The first which is the original from \cite{C-S}, is
quite helpful in understanding Definition~\ref{D2.4}. Let $e_j$ be the
minimal projection in $D_i$ which is 1 in the $j$'th row. Let
\[
x_{i_1\ldots i_k}=\cdots\otimes 1\otimes e_{i_1}\otimes\cdots\otimes
e_{i_k}\otimes 1\in D_{1k}\;.
\]
Then $(x_{i_1,\ldots,i_k})\in S_k$, and
\[
x_{i_j}^i=\cdots\otimes1\otimes e_{i_j}\otimes1\otimes
    \cdots\in D_j\;.
\]
Thus
\begin{eqnarray*}
H(M_1,\alpha(M_1),\ldots,\alpha^{k-1}(M_1)) &=&
     H(M_1,\ldots,M_k) \\
&\geq& \sum_{i_1\ldots i_k}\eta\tau(x_{i_1\ldots i_k})
   -\sum_{j=1}^k \sum_{i_j=1}^n \tau\eta(E_{M_j}x_{i_j}^j) \\
&=& n^k \eta(n^{-k})-\sum_j \sum_{i_j} \tau\eta(e_{i_j}) \\
&=& k\log n-0\;,
\end{eqnarray*}
so that $H(M_1,\alpha)\geq\log n$.

To prove the opposite inequality we use that $(M_1\cup
M_2\cup\cdots\cup M_k)''$ is a factor of type I$_{n^k}$, hence has
entropy
$k\log n$. The rest of the proof consists of an application of the
Kolmogoroff-Sinai theorem to the sequence $(P_j)$ together with an
application of property $(E)$.

The other proof is quicker. Fix $q\in\N$. Then $A_q=D_{-q,q}$ is a masa
in $P_q$. If $k\in\N$ let $A=\bigvee\limits_{j=0}^{k-1}\alpha^j(A_q)$.
Then $A=D_{-q,q+k-1}$ is a masa in
$\bigvee\limits_0^{k-1}\alpha^j(P_q)$ such that
$A\cap\alpha^j(P_q)=\alpha^j(A_q)$ is a masa in
$\alpha^j(P_q)$. Thus by the easy part of Theorem~\ref{T2.6},
\[
H(P_q,\alpha)=\lim_{k\to\infty}\frac{2q+k-1}{k}\log n=\log n\;,
\]
so that by the Kolmogoroff-Sinai theorem, $H(\alpha)=\log n$.
\end{remark}

\paragraph{2.14  Bernoulli shifts} \
The above arguments can be extended to noncommutative Bernoulli shifts
of the hyperfinite II$_1$-factor $R$. Let $h\in M_0^+$ with $\Tr(h)=1$,
$\Tr$ denoting the usual trace on $M_n(\C)$, with eigenvalues
$h_1,\ldots,h_n$. Let $\varphi_0$ be the state $\varphi_0(x)=\Tr(hx)$
for $x\in M_0$.

Let $\varphi_i=\varphi_0$ on $M_i$ and
$\varphi=\bigotimes\limits_{i\in\Z} \varphi_i$ denote the corresponding
product state on $B=\bigotimes\limits_{i\in\Z} M_i$. In the
GNS-representation of $B$ due to $\varphi$ the centralizer $R$ of the
weak closure is the hyperfinite II$_1$-factor and contains the algebras
$A_q$ above. Since $\alpha$ is $\varphi$-invariant, the extension of
$\alpha$ to the GNS-representation restricts to an automorphism of $R$,
which we call a Bernoulli shift. A slight extension of the argument
from 2.13 shows that
\[
H(\alpha)=\sum_{i=1}^n \eta(h_i)=S(\varphi_0)\;,
\]
where $S(\varphi_0)$ is the entropy of the state $\varphi_0$ on
$M_n(\C)$.
\hfill
$\Box$
\bigskip

There is another natural definition of Bernoulli shift on $R$. Let $T$
be a (classical) Bernoulli shift on a probability space $(X,{\mathcal
B},\mu)$. Then $R=L^\infty(X,{\mathcal B},\mu)\times_T\Z$, and $T$
extends to an inner automorphism $\Ad u_T$ on $R$. Both the von Neumann
algebra generated by the $A_q$'s above and $L^\infty(X,{\mathcal
B},\mu)$ are Cartan subalgebras of $R$, i.e. they are masas whose
normalizers generate $R$, hence they are conjugate by \cite{CFW}. Thus
we have one outer and one inner automorphism on $R$ which act as the
same Bernoulli shift on a Cartan subalgebra and have the same entropy
(see Theorem~7.1 below).

Other extensions of classical shift automorphisms have been studied by
Besson \cite{B} for Markov shifts and Quasthoff \cite{Q}. In all
examples there is a masa like $A$ in 2.19 and the entropy is the same
as the classical counterpart.

\section{Entropy in C*-algebras}

After the appearance of \cite{C-S} an obvious problem was to extend the
definition from the tracial case to that of general states. It took 10
years before Connes \cite{Co} saw what had to be done. If one looks at
the rewritten form of $H(N_1,\ldots,N_k)$ in equation 2.5 and notes
that states are of the form $\varphi(x)=\tau(hx)$, $h\in
L^1(M,\tau)^+$, it is obvious what to do. Let $\varphi$ be a normal
state of a von Neumann algebra $M$. Modify Notation (2.2) as follows:
For
$k\in\N$ put
\begin{eqnarray*}
&&
S_{k,\varphi}=\{
\varphi_{i_1\ldots i_k}\in M_\ast^+\;,\quad i_j\in\N\;,\quad
     \varphi_{i_1\ldots i_k}=0\quad \mbox{except for} \\
&&\qquad\qquad \mbox{a finite number of indices},\quad
     \sum_{i_1\ldots i_k}\varphi_{i_1\ldots i_k}=\varphi\}\;.
\end{eqnarray*}
Let
\[
\varphi_{i_j}^j=\sum_{i_1\ldots i_{j-1}i_{j+1}\ldots i_k}
       \varphi_{i_1\ldots i_k}\;.
\]
If $N_1,\ldots,N_k\subset M$ are finite dimensional von Neumann
subalgebras we let \cite{Co}
\begin{equation}\label{e3.1}
\qquad
H_\varphi(N_1,\ldots,N_k)=
\sup_{(\varphi_{{i_j}\ldots i_k})\in S_{k,\varphi}}
\Big\{ \sum_{i_1\ldots i_k}\eta(\varphi_{i_1\ldots i_k}(1))
+\sum_{j=1}^k \sum_{i_j} S(\varphi_{i_j}^j|_{N_j},
\varphi|_{N_j})\Big\}\;.
\end{equation}
This definition even makes sense for C*-algebras, because, as we
pointed out earlier, Pusz, Woronowicz and Kosaki extended the
definition of relative entropy to C*-algebras. Since C*-algebras may
have no finite dimensional C*-subalgebras except the scalars, the
definition above would only be useful for AF-algebras and their like.
Connes, together with Narnhofer and Thirring \cite{CNT} circumvented
the problem by replacing the algebras $N_j$ by completely positive maps
$\gamma_j$ from finite dimensional algebras into the C*-algebra. The
definition is as follows.

Let $A$ be a unital C*-algebra with a state $\varphi$. Let
$N_1,\ldots,N_k$ be finite dimensional C*-algebras and $\gamma_j:N_j\to
A$ a unital completely positive map, $j=1,\ldots,k$. Let $B$ be a
finite dimensional abelian C*-algebra and $P:A\to B$ a unital positive
linear map such that there is a state $\mu$ on $B$ with $\mu\circ
P=\varphi$. Let $p_1,\ldots,p_r$ be the minimal projections in $B$.
Then there are states $\hat{\varphi}_1,\ldots,\hat{\varphi}_r$ on $A$
such that
\begin{equation}\label{e3.2}
P(x)=\sum_{i=1}^r\hat{\varphi}_i(x)p_i\;,
\end{equation}
and
\begin{equation}\label{e3.3}
\varphi=\sum_{i=1}^r \mu(p_i)\hat{\varphi}_i\;,
\end{equation}
is $\varphi$ written as a convex sum of states. Put
\[
\varepsilon_\mu(P)=\sum \mu(p_i)S(\hat{\varphi}_i,\varphi)\;.
\]
where $S(\hat{\varphi}_i,\varphi)$ is the relative entropy. Let the
{\em entropy defect\/} be
\[
s_\mu(P)=S(\mu)-\varepsilon_\mu(P)\;,
\]
where $S(\mu)=\sum\limits_{i=1}^r \eta(\mu(p_i))$ is the entropy of
$\mu$.

Suppose $B_1,\ldots,B_k$ are C*-subalgebras of $B$ and $E_j:B\to B_j$
the $\mu$-invariant conditional expectaion. Then the quadruple
$(B,E_j,P,\mu)$ is called an {\em abelian model\/} for
$(A,\varphi,\gamma_1,\ldots,\gamma_k)$, and its entropy is defined to be
\begin{equation}\label{e3.4}
S\Big(\mu\Big|_{\bigvee\limits_{j=1}^k B_j}\Big)
-\sum_{j=1}^k s_\mu(P_j)\;,
\end{equation}
where $P_j=E_j\circ P\circ\gamma_j:N_j\to B_j$, and the definition of
$s_\mu(P_j)$ is the same as for $P$ above, where we replace $\mu$ by
$\mu|_{B_j}$, $\varphi$ by $\varphi\circ \gamma_j$.

\begin{definition}\label{D3.5}
{\rm \cite{CNT}} \
$H_\varphi(\gamma_1,\ldots,\gamma_k)=\sup$ of (\ref{e3.4}) over all
abelian models.
\end{definition}

In the special case when $N_1,\ldots,N_k\subset A$ and $\gamma_j:N_j\to
A$ is the inclusion map let $(B,E_j,P,\mu)$ be an abelian model. We may
assume $B=\bigvee\limits_{j=1}^k B_j$. Let $\{p_{i_j}\}$ be the set of
minimal projections in $B_j$. Then $\{p_{(i)}=p_{i_1\ldots i_k}=
p_{i_1}\ldots p_{i_k}\}$ is the set of minimal projections in $B$. If
we use the abbreviation following Definition~\ref{D2.4}, equations
(\ref{e3.2}) and (\ref{e3.3}) can be written
\begin{eqnarray*}
P(x) &=& \sum_{(i)}\hat{\varphi}_{(i)}(x)p_{(i)}\;, \\
\varphi &=& \sum_{(i)}\mu(p_{(i)})\hat{\varphi}_{(i)}\;.
\end{eqnarray*}
Let $\varphi_{(i)}(x)=\mu(p_{(i)})\hat{\varphi}_i(x)$, so $\varphi=\sum
\varphi_{(i)}$. We have
\[
E_j(p_{(i)})=\mu(p_{(i)})\mu(p_{i_j}^j)^{-1}p_{i_j}^j\;,
\]
where $p_{i_j}=p_{i_j}^j$ in the notation of (2.2). Hence if
$x\in N_j$ we have
\[
P_j(x)=E_j\circ P(x)=\sum_{(i)}
\frac{\varphi_{(i)}^j(x)}{\mu(p_{i_j})}\mu(p_{(i)})p_{i_j}=
\sum_{i_j}\frac{\varphi_{i_j}^j(x)}{\mu(p_{i_j})}p_{i_j}\;.
\]
Thus
\[
\varepsilon_\mu(P_j)=\sum_{i_j}\mu(p_{i_j})S\Big(
\frac{\varphi_{i_j}^j}{\mu(p_{i_j})}\Big|_{N_j},
\varphi|_{N_j}\Big)\;.
\]
>From the identity
\[
\lambda S(\rho,\psi)=\eta(\lambda)+S(\lambda\rho,\psi)\;,
\]
as is easily shown in the finite dimensional case, we have
\begin{eqnarray*}
\varepsilon_\mu(P_j) &=&
\sum_{i_j}\{ \eta(\mu(p_{i_j}))+S(\varphi_{i_j}^j|_{N_j},
\varphi|_{N_j})\} \\
&=& S(\mu|_{B_j})+\sum_{i_j} S(\varphi_{i_j}^j|_{N_j},
\varphi|_{N_j})\;.
\end{eqnarray*}
Note that $\varphi_{(i)}(1)=\mu(p_{(i)})$. Thus we find that
the entropy of the abelian model $(B,E,P,\mu)$ is
\[
S(\mu|_B)-\sum_{j=1}^k \{ S(\mu|_{B_j})-
\varepsilon_\mu(P_j)\}
=
\sum_{(i)}\eta(\varphi_{(i)}(1))
+\sum_{j=1}^k \sum_{i_j} S(\varphi_{i_j}^j|_{N_j},
\varphi|_{N_j})\;,
\]
which is the same as the expression in (\ref{e3.1}).

We note that if we are given $(\varphi_{(i)})\in S_{k,\varphi}$ and
$N_1,\ldots,N_k\subset M$ it is not hard to construct an abelian model
like $(B,E,P,\mu)$ above, so that (\ref{e3.1}) defines
$H_\varphi(\gamma_1,\ldots,\gamma_k)$ when $\gamma_j:N_j\to M$ is the
inclusion map. Thus $H_\varphi(\gamma_1,\ldots,\gamma_k)$ is a direct
generalization of $H_\tau(N_1,\ldots,N_k)$ defined in
Definition~\ref{D2.4}. One can show similar properties to (A)--(E) in
Section~2, see
\cite{CNT}, hence we can define the entropy of an automorphism

\begin{definition} \label{D3.6}
{\rm \cite{CNT}} \
Let $A$ be a unital C*-algebra, $\varphi$ a state and $\alpha$ a
$\varphi$-invariant automorphism of $A$. Let $C$ be a finite
dimensional C*-algebra and $\gamma:C\to A$ a unital completely positive
map. Then
\[
h_{\varphi,\alpha}(\gamma)=\lim_{k\to\infty}\frac{1}{k}
H_\varphi(\gamma,\alpha\circ \gamma,\ldots,\alpha^{k-1}\circ \gamma)
\]
exists. We define the {\em entropy of\/} $\alpha$ with respect to
$\varphi$ to be
\[
h_\varphi(\alpha)=\sup_{(C,\gamma)}h_{\varphi,\alpha}(\gamma)\;,
\]
where the $\sup$ is taken over all pairs $(C,\gamma)$.
\end{definition}

The Kolmogoroff-Sinai Theorem, cf. Theorem~2.11, now takes the
following form.

\begin{theorem} \label{T3.7}
{\rm \cite{CNT}} \
Let $A,\varphi,\alpha$ be as above. Suppose $(\tau_n)$ is a sequence of
unital completely positive maps $\tau_n:A_n\to A$ from finite
dimensional C*-algebras $A_n$ into $A$ such that there exist unital
completely positive maps $\sigma_n:A\to A_n$ for which $\tau_n\circ
\sigma_n\to\id_A$ in the pointwise norm topology. Then
\[
\lim_{n\to\infty}h_{\varphi,\alpha}(\tau_n)=h_\varphi(\alpha)\;.
\]
In particular if $A=\overline{\bigcup\limits_{n=1}^\infty A_n}$ is an
AF-algebra, and we idenitfy $A_n$ with its inclusion map $A_n\to A$, we
have,
\[
h_\varphi(\alpha)=\lim_{n\to\infty} h_{\varphi,\alpha}(A_n)\;.
\]
\end{theorem}

Theorem~\ref{T3.7} is applicable if $A$ is nuclear. In that case we
have \cite{CNT}
\begin{itemize}
\item[(i)]
$h_\varphi(\alpha)=h_\varphi(\gamma\alpha\gamma^{-1})$ if
$\gamma\in\Aut A$.

\item[(ii)]
$h_\varphi(\alpha^n)=|n|(h_\varphi(\alpha)$, $n\in\Z$.

\item[(iii)]
If $\varphi_1$ and $\varphi_2$ are $\alpha$-invariant and
$\lambda\in[0,1]$,
$h_{\lambda\varphi_1+(1-\lambda)\varphi_2}(\alpha)\geq\lambda
h_{\varphi_1}(\alpha)+(1-\lambda)h_{\varphi_2}(\alpha)$.
\end{itemize}

Definition~\ref{D3.6} also makes sense for normal states of von Neumann
algebras. If $(A,\varphi,\alpha)$ is as above, and
$(\pi_\varphi,H_\varphi,\xi_\varphi)$ is the GNS-representation of
$\varphi$, and $\bar{\alpha}$ the extension of $\alpha$ to
$M=\pi_\varphi(A)''$, then \cite{CNT},
\begin{equation} \label{e3.8}
h_\varphi(\alpha)=h_{\omega_{\xi_\varphi}}(\bar{\alpha})\;.
\end{equation}
Thus we can freely move back and forth between $A$ and $M$ in our
computations of entropy.

If $(A_i,\varphi_i,\alpha_i)$, $i=1,2$ are C*-dynamical systems we have
as in the tracial case (2.6),
\begin{equation} \label{e3.9}
h_{\varphi_1\otimes\varphi_2}(\alpha_1\otimes\alpha_2)\geq
h_{\varphi_1}(\alpha_1)+h_{\varphi_2}(\alpha_2)\;,
\end{equation}
because there are many more choices of abelian models to compute the
left side of (\ref{e3.9}) than the right.

In nicer cases a useful criterion for computing entropy on von Neumann
algebras is to consider the restriction of the automorphism to the
centralizer of the state.

\begin{proposition} \label{P3.10}
{\rm \cite{CNT}} \
Let $M$ be a von Neumann algebra and $\varphi$ a normal state. Let
$N_1,\ldots,N_k$ be finite dimensional von Neumann subalgebras of $M$.
Suppose they contain abelian subalgebras $A_j\subset N_j\cap
M_\varphi$, where $M_\varphi$ is the centralizer in $M$, such that the
$A_j$ pairwise commute and $A=\bigvee\limits_{j=1}^k A_j$ is a masa in
$N=\bigvee_{j=1}^k N_j$. Then
\[
H_\varphi(N_1,\ldots,N_k)=S(\varphi|_N)\;.
\]
\end{proposition}

A good illustration of Proposition~\ref{P3.10} is the case of Bernoulli
shifts as described in Theorem~2.14. In the notation of 2.14 the
state $\varphi=\bigotimes\limits_{i\in\Z}\varphi_i$ on
$B=\bigotimes\limits_{i\in\Z}M_i$ satisfies the conditions of the
proposition, so with $\alpha$ the shift,
$h_\varphi(\alpha)=H_{\varphi|_{M_\varphi}}(\alpha|_{M_\varphi})
=S(\varphi_0)$.

\section{Bogoliubov automorphisms}

The main examples for which the C*-algebra entropy have
been computed, are those of quasifree states of the CAR- and
CCR-algebras and invariant Bogoliubov (or quasifree) automorphisms.
The computations and results are quite similar, so for simplicity we
restrict attention to the CAR-algebra. Let us recall the definitions.

Let $H$ be a complex Hilbert space. The CAR-algebra ${\cal A}
(H)$ over $H$ is a C*-algebra with the property that
there is a linear map $f\to a(f)$ of $H$ into ${\cal A}(H)$
whose range generates ${\cal A}(H)$ as a C*-algebra and
satisfies the canonical anticommutation relations
\begin{eqnarray*}
&&a(f)a(g)^\ast+a(g)^\ast a(f)=(f,g)1,\qquad
          f,g\!\in\!H\,, \\
&&a(f)a(g)+a(g)a(f)=0\,,
\end{eqnarray*}
where $(\cdot,\cdot)$ is the inner product on $H$ and 1 the
unit of ${\cal A}(H)$. If $0\leq A\leq1$ is an operator on $H$,
then the quasifree state $\omega_A$ on ${\cal A}(H)$ is defined
by its values on products of the form $a(f_n)^\ast\ldots
a(f_1)^\ast a(g_1)\ldots a(g_m)$ given by
\[
\omega_A(a(f_n)^\ast\ldots a(f_1)^\ast a(g_1)\ldots
a(g_m))=\delta_{nm}\det((Ag_i,f_j))\,.
\]
If $U$ is a unitary operator on $H$ then $U$ defines an
automorphism $\alpha_U$ on ${\cal A}(H)$, called a {\em
Bogoliubov automorphism}, determined by
\[
\alpha_U(a(f))=a(Uf)\,.
\]
If $U$ and $A$ commute it is an easy consequence of the above
definition of $\omega_A$ that $\alpha_U$ is
$\omega_A$-invariant. We shall in this section state a formula for the
entropy
$h_{\omega_A}(\alpha_U)$ and indicate the ideas in the computation
when $A$ is a scalar operator.

Connes suggested to the author that if
$\omega_A$ is the trace $\tau$ the answer should be
\begin{equation} \label{e4.1}
h_\tau(\alpha_U)=\frac{\log 2}{2\pi}\int_0^{2\pi}
m(U)(\theta)d\theta\,,
\end{equation}
where $m(U)$ is the multiplicity function of the absolutely
continuous part $U_a$ of $U$. Then Voiculescu and St\o{}rmer \cite{SV}
showed this and more by solving the problem when $A$ has pure point
spectrum. Later on Narnhofer and Thirring \cite{N-T1} and Park and
Shin \cite{P-S} independently extended the result to more general
$A$. Finally Neshveyev \cite{N} settled the problem completely in the
general case. He and Golodets \cite{G-N1}, and before them Bezuglyi and
Golodets \cite{B-G} considered more general group actions than $\Z$.

If $A$ has pure point spectrum there is an orthonormal basis
$(f_n)$ of $H$ such that $Af_n=\lambda_n f_n$,  $n\!\in\!\Bbb
N$, $0\!\leq\!\lambda_n\!\leq\!1$. Define recursively
operators
\begin{eqnarray*}
&&V_0=1,\qquad V_n=\prod_{i=1}^n\Big(1-2a(f_i)^\ast
     a(f_i)\Big),\qquad e_{11}^{(n)}=a(f_n)a(f_n)^\ast \\
&&e_{12}^{(n)}=a(f_n)V_{n-1},\quad e_{21}^{(n)}=
V_{n-1}a(f_n)^\ast,\quad e_{22}^{(n)}=a(f_n)^\ast a(f_n)\,.
\end{eqnarray*}
Then the $e_{ij}^{(n)}$, $i,j\!=\!1,2$ form a complete set of
$2\times2$ matrix units generating a I$_2$-factor $M_2({\Bbb
C})_n$, and for distinct $n$ and $m$ $e_{ij}^{(n)}$ and
$e_{k\ell}^{(m)}$ commute. Thus ${\cal
A}(H)\simeq\bigotimes\limits_1^\infty M_2({\Bbb C})_n$, and
$\omega_A$ is a product state
$\omega_A=\bigotimes\limits_1^\infty\omega_{\lambda_n}^0$ with
respect to this factorization, where $\omega_\lambda^0$ is the
state on $M_2({\Bbb C})$ given by
\[
\omega_\lambda^0\left(\begin{pmatrix}
a & b \\ c & a\end{pmatrix}\right)=
(1-\lambda)a+\lambda d\,.
\]
In case $A=\lambda1$ we write $\omega_\lambda$ for $\omega_A$.
Then $\alpha_U$ is $\omega_\lambda$-invariant for all $U$. We
consider the entropy $h_{\omega_\lambda}(\alpha_U)$.

Each unitary $U$ is a direct sum $U=U_a\oplus U_s$, where $U_a$
has spectral measure absolutely continuous with respect to
Lebesgue measure $d\theta$ on the circle, while $U_s$ has
spectral measure singular with respect to $d\theta$. We shall
as above denote by $m(U)$ the multiplicity function of $U_a$.
The idea is now to approximate the case  when
\[
U=U_s\oplus U_1\oplus\cdots\oplus U_n\,,
\]
where each $U_i$ acts on a Hilbert space $H_i$, $i\!=\!1,
\ldots,n$, and $U_i$ is unitarily equivalent to $V^{p_i}$,
where $V$ is a bilateral shift. Let us for simplicity ignore
the complications due to the grading of ${\cal A}(H)$ as a
direct sum of its even and odd parts. Then
\[
\alpha_U=\alpha_{U_s}\otimes\alpha_{U_1}\otimes
\cdots\otimes\alpha_{U_n}
\]
and
\[
\omega_\lambda=\omega_\lambda|{\cal
A}(H_s)\otimes\omega_\lambda|{\cal
A}(H_1)\otimes\cdots\otimes\omega_\lambda|{\cal A}(H_n)\,.
\]
Thus we could hope that
\begin{equation} \label{e4.2}
h_{\omega_\lambda}(\alpha_U)=h_{\omega_\lambda|{\cal A}(H_s)}
(\alpha_{U_s})+\sum_{i=1}^n h_{\omega_\lambda|
    {\cal A}(H_i)}(\alpha_{U_i}),
\end{equation}
and thus restrict attention to the case when $U$ is singular
or a power of a bilateral shift. We do have problems because of
(\ref{e3.9}), but life turns out nicely because we can as with
the shift in 2.13 restrict attention to the diagonal, and
the diagonal is contained in the even CAR-algebra, where the
tensor product formulas above hold. First we take care of the
singular part $U_s$.

\begin{lemma}\label{L4.1}
If $U$ has spectral measure singular with respect to the
Lebesgue measure, and $\alpha_U$ is $\varphi$-invariant for a
state $\varphi$, then $h_\varphi(\alpha_U)=0$.
\end{lemma}

Thus in (\ref{e4.2}) we can forget about $U_s$. If $U=V^p$ with $V$ a
bilateral shift and $p\!\in\!\Bbb Z$, then
\[
h_{\omega_\lambda}(\alpha_U)=h_{\omega_\lambda}((\alpha_V)^p)=
|p|h_{\omega_\lambda}(\alpha_V)\,.
\]
If we write ${\cal A}(H)=\bigotimes\limits_{n=-\infty}^\infty
M_2({\Bbb C})_n$, then on the diagonal $\alpha_V$ is the
shift, so like in 2.14, see also Proposition~\ref{P3.10}, we get
\[ h_{\omega_\lambda}(\alpha_V)=\eta(\lambda)+
     \eta(1-\lambda)\,.
\]
Now $|p|$ is the multiplicity $m(U)$ of $U$, and since
$\frac{1}{2\pi}d\theta$ is the normalized Haar measure on the
circle, it is not surprising that we have

\begin{theorem} \label{T4.4}
Let $U$ be a unitary operator on $H$ and $\lambda\in[0,1]$.
Then $h_{\omega_\lambda}(\alpha_U)=\frac{1}{2\pi}\Big(
\eta(\lambda)+\eta(1-\lambda)\Big) \int_0^{2\pi}
m(U)(\theta)d\theta$.
\end{theorem}

Note that if $\lambda=1/2$, $\omega_\lambda=\tau$, so we get
formula (\ref{e4.1}). For more general $A$ we use direct integral
theory with respect to the von Neumann algebra generated by
$U_a$. If $A$ commutes with $U$, $A=A_a\oplus A_s$, where
$A_a=\int\limits_0^{2\pi\oplus} A(\theta)d\theta$,
$H=\int\limits^\oplus H_\theta d\theta$, and $H_\theta=0$ if
$m(U)(\theta)=0$, and $A(\theta)\subset B(H_\theta)$.

We can now state the main theorem in this section, see \cite{SV},
\cite{P-S}, \cite{N-T1} and \cite{N}.

\begin{theorem} \label{T4.5}
Let $0\leq A\leq1$ and $U$ be a unitary operator
commuting with $A$. Then
\[
h_{\omega_A}(\alpha_U)=\frac{1}{2\pi}
\int\limits_0^{2\pi}
Tr(\eta(A(\theta))+\eta(1-A(\theta)))d\theta\,.
\]
Furhtermore, \cite{N}, $h_{\omega_A}(\alpha_U)<\infty$ if and only if
$A(\theta)$ has pure point spectrum for almost all $\theta\in[0,2\pi)$.
\end{theorem}

Entropy is far from a complete conjugacy invariant for Bogoliubov
automorphisms. It has been shown by Golodets and Neshveyev \cite{G-N3}
in the case of the CCR-algebra, that there exists a quasifree state
$\omega$ and a one-parameter family $\alpha_\theta$,
$\theta\in[0,2\pi)$, of $\omega$-invariant Bogoliubov automorphisms
with the same positive entropy $h_\omega(\alpha_\theta)$ and such that
if $M$ is the weak closure of the CCR-algebra in the GNS-representation
of $\omega$, then the W*-dynamical systems $(M,\omega,\alpha_\theta)$
are pairwise nonconjugate. In this case $M$ is a factor of type III$_1$
and the centralizer $M_\omega$ of $\omega$ in $M$ is the scalars. Thus
the situation is quite different from that encountered in
Proposition~\ref{P3.10}. Such an example had previously been found by
Connes \cite{Co}.

\section{The entropy of Sauvageot and Thouvenot}

Sauvageot and Thouvenot \cite{S-T} have given an alternative definition
of entropy, which is close to that of \cite{CNT}, but which has
technical advantages in some cases. Again we look at all possible ways
a state can be written as convex combinations of other states.

Let $(A,\varphi,\alpha)$ be a unital C*-dynamical system, and let
$(C,\mu,\beta)$ be an abelian C*-dynamical system. A {\em stationary
coupling\/} of these two systems is an $\alpha\otimes\beta$-invariant
state $\lambda$ on $A\otimes C$ such that $\lambda|_A=\varphi$,
$\lambda|_C=\mu$. If $P$ is a finite dimensional C*-subalgebra of $C$
with atoms $p_1,\ldots,p_r$ let
\[
\varphi_i(x)=\mu(p_i)^{-1}\lambda(x\otimes p_i)\;,
\]
whenever $\mu(p_i)\not=0$, as we may assume. Then
\[
\varphi=\sum_{i=1}^r \mu(p_i)\varphi_i
\]
is $\varphi$ written as a convex sum of states. Let
\[
P^-=\bigvee_{i=1}^\infty \beta^i(P)\;;
\]
then from the classical theory
\[
H_\mu(P,\beta)=\lim_{n\to\infty}\frac{1}{n}H_\mu
\Big( \bigvee_{i=0}^{n-1}\beta^i(P)\Big)=
H_\mu(P|P^-)\;.
\]

\begin{definition} \label{D5.1}
The {\em Sauvageot-Thouvenot entropy\/} $h'_\varphi(\alpha)$ of the
system $(A,\varphi,\alpha)$ is the supremum over all abelian systems
$(C,\mu,\beta)$ as above of the quantities
\[
H_\mu(P|P^-)-H_\mu(P)+
\sum_{i=1}^r \mu(p_i)S(\varphi_i,\varphi)\;.
\]
\end{definition}

Note that if $\lambda=\varphi\otimes\mu$ then $\varphi_i=\varphi$,
hence $S(\varphi_i,\varphi)=0$, and $H(P|P^-)-H_\mu(P)\leq0$. Thus if
$\lambda=\varphi\otimes\mu$ is the only stationary coupling then
$h'_\varphi(\alpha)=0$, a fact we shall need later.

\begin{theorem} \label{T5.2}
{\rm \cite{S-T}} \
If $A,\varphi,\alpha)$ is a C*-dynamical system with $A$ nuclear or a
W*-dynamical system with $A$ injective then $h'_\varphi(\alpha)$
equals the CNT-entropy $h_\varphi(\alpha)$.
\end{theorem}

In order to understand the relation between $h_\varphi$ and
$h'_\varphi$ better, let us prove the inequality
$h'_\varphi(\theta)\leq h_\varphi(\alpha)$ in some detail. The
opposite inequality holds in general.

Let $(C,\mu,\beta)$ be an abelian system and $\lambda$ a stationary
coupling. Let $P$ be a finite dimensional subalgebra of $C$ with atoms
$p_1,\ldots,p_r$. Let $m\in\N$. We shall show
\begin{equation}\label{e5.3}
h_\varphi(\alpha)\geq H_\mu(P|P^-)
-H_\mu(P)+\sum_{i=1}^r \mu(p_i)S(\varphi_i,\varphi)\;,
\end{equation}
where $\varphi_i(x)=\mu(p_i)^{-1}\lambda(x\otimes p_i)$.

$\varphi$ has the decomposition
\[
\varphi=\sum_{i_1,\ldots,i_m}\varphi_{i_1\ldots i_m}\;,
\]
where
\[
\varphi_{i_1\ldots i_m}(x)=
\lambda(x\otimes p_{i_1}\beta(p_{i_2})\ldots \beta^{m-1}(p_{i_m}))\;.
\]
We thus have a completely positive map $\rho: A\to C$ defined by
\[
\rho(x)=\sum
\frac{\varphi_{i_1\ldots i_m}(x)}{\mu(p_{i_1\ldots i_m})}
p_{i_1\ldots i_m}\;,
\]
where $p_{i_1\ldots i_m}=p_{i_1}\beta(p_{i_2})\ldots
\beta^{m-1}(p_{i_m})$ are the atoms in
$B=\bigvee\limits_0^{m-1}\beta^j(P)$. The $\mu$-invariant conditional
expectation $E_j:B\to\beta^j(P)$ is determined by
\[
E_j(p_{i_1\ldots i_m})=
\frac{\mu(p_{i_1\ldots i_m})}{\mu(\beta^j(p_{i_j}))}\beta^j(p_{i_j})\;.
\]

Let $\gamma$ be a unital completely positive map from a finite
dimensional C*-algebra into $A$. We found above an abelian model
for $(A,\varphi,\gamma,\alpha\circ\gamma,\ldots,\alpha^{m-1}
\circ\gamma)$, namely $(B,E_j,\rho,\mu)$. We find
\begin{eqnarray*}
E_j\circ\rho(x) &=& \sum_{i_1\ldots i_m}
   \frac{\varphi_{i_1\ldots i_m}(x)}{\mu(\beta^j(p_{i_j}))}
     \beta^j(p_{i_j}) \\
&=& \sum_{i_j}
    \frac{\lambda(x\otimes\beta^j(p_{i_j}))}{\mu(\beta^j(p_{i_j}))}
    \beta^j(p_{i_j}) \\
&=& \sum_i
    \frac{\lambda(\alpha^{-j}(x)\otimes p_i)}{\mu(p_i)}
    \beta^j(p_i) \\
&=& \sum_i \varphi_i(\alpha^{-j}(x))\beta^j(p_i)
\end{eqnarray*}
by invariance of $\mu$ and $\lambda$ with respect to $\beta$ and
$\alpha\otimes\beta$ respectively. From Definition~\ref{D3.5} we obtain
\begin{eqnarray*}
\lefteqn{
H_\varphi(\gamma,\alpha\circ\gamma,\ldots,\alpha^{m-1}\circ\gamma)
    \geq S(\mu|B)-\sum_{j=0}^{m-1}S(\mu|\beta^j(P))} \\
&&\qquad + \sum_{j=0}^{m-1} \sum_{i=1}^r
     \mu(\beta^j(p_i))S(\varphi_i\circ\alpha^{-j}
     \circ\alpha^j\circ\gamma,
    \varphi\circ\alpha^j\circ\gamma) \\
&&=S(\mu|B)-mS(\mu|P) + m\sum_{i=1}^r \mu(p_i)
    S(\varphi_i\circ\gamma,\varphi\circ\gamma)\;.
\end{eqnarray*}
Thus
\[
\frac{1}{m}H_\varphi(\gamma,\alpha\circ\gamma,\ldots,\alpha^{m-1}
\circ\gamma)\geq
H_\mu(P|P^-)-H_\mu(P)
+\sum_i \mu(p_i)S(\varphi_i\circ\gamma,\varphi\circ\gamma)\;.
\]
If $A$ is an injective von Neumann algebra or a nuclear C*-algebra we
can find a net of maps $(\gamma_\omega)_\omega$ such that
$S(\varphi_i\circ\gamma_\omega,\varphi\circ\gamma_\omega)\to
S(\varphi_i,\varphi)$, see \cite[5.29 and 5.30]{O-P}. Thus (\ref{e5.3})
holds, and $h_\varphi(\alpha)\geq h'_\varphi(\alpha)$.
\hfill
$\Box$

\bigskip

If $\alpha$ is an automorphism of a C*-algebra $A$, and $B$ is a
C*-subalgebra of $A$ such that $\alpha(B)=B$, then every
$\alpha$-invariant state on $B$ has an $\alpha$-invariant extension to
$A$. If we apply this to a stationary coupling $\lambda$ on $N\otimes
C$ we can extend $\lambda$ to $A\otimes C$ and prove the following
result.

\begin{proposition}\label{P5.3}
{\rm \cite{N-S1}} \
Let $A$ be a unital C*-algebra and $\alpha$ an automorphism of $A$. Let
$B$ be an $\alpha$-invariant C*-subalgebra of $A$ and $\psi$ an
$\alpha$-invariant state of $B$. Then for each $\varepsilon>0$ there
exists an $\alpha$-invariant state $\varphi$ of $A$ such that
$\varphi|_B=\psi$ and
$h_\varphi'(\alpha)>h_\psi'(\alpha|_B)-\varepsilon$.
\end{proposition}

It is well-known from subfactor theory that properties of a subfactor
of finite index often are kept by the larger factor. The following is
an analogous result for entropy.

\begin{proposition}\label{P5.4}
{\rm \cite{N-S1}} \
Let $(A,\varphi,\alpha)$ be a unital C*-dynamical system. Let $B\subset
A$ be an $\alpha$-invariant C*-subalgebra with $1\in B$. Suppose there
is a conditional expectation $E:A\to B$ such that
$E\circ\alpha=\alpha\circ E$, $\varphi\circ E=\varphi$, and $E(x)\geq
cx$ for all $x\in A^+$ for some real number $c>0$. Then
\[
h'_\varphi(\alpha)=h'_\varphi(\alpha|_B)\;.
\]
\end{proposition}

\begin{proof}
The inequality $h'_\alpha(\alpha)\geq h'_\varphi(\alpha|_B)$ follows
by monotonicity. To prove the opposite inequality we consider $A$ in its
GNS-representation with respect to $\varphi$, so we may assume $A$ and
$B$ are von Neumann algebras and $\varphi$ and $E$ normal. Let
$(C,\mu,\beta)$ be as before. Then $\varphi_i\circ E\geq c\varphi_i$,
hence by (\ref{e2.7})
$S(\varphi_i,\varphi_i\circ E)\leq S(\varphi_i,c\varphi_i)=-\log c$.
Thus by (\ref{e2.6})
\begin{eqnarray*}
\sum\mu(p_i)S(\varphi_i,\varphi) &=& \sum\mu(p_i)
    (S(\varphi_i|_B,\varphi|_B)
     +S(\varphi_i,\varphi_i\circ E)) \\
&\leq& \sum\mu(p_i)S(\varphi_i|_B),\varphi|_B)-\log c
\end{eqnarray*}
It follows from Definition~\ref{D5.1} that $h'_\varphi(\alpha)\leq
h'_\varphi(\alpha|_B)-\log c$. But this inequality must hold for
$\alpha^m$, $m\in\N$, as well. Hence
\[
h'_\varphi(\alpha)=\frac{1}{m} h'_\varphi(\alpha^m)\leq
\frac{1}{m}h'_\varphi(\alpha^m|_B)-
\frac{1}{m}\log c=
h'_\varphi(\alpha|_B)-
\frac{1}{m}\log c\;,
\]
and $h'_\varphi(\alpha)\leq h'_\varphi(\alpha|_B)$.
\hfill
$\Box$
\end{proof}
\bigskip

>From Theorem~\ref{T5.2} it follows that if $A$ and $B$ are nuclear
C*-algebras or injective von Neumann algebras then
\[
h_\varphi(\alpha)=h_\varphi(\alpha|_B)\;.
\]

\begin{theorem}\label{T5.5}
{\rm \cite{N-S1}} \
If $(M,\varphi,\alpha)$ is a W*-dynamical system with $M$ of type I
and $Z$ is the center of $M$, then
$h_\varphi(\alpha)=h_\varphi(\alpha|_Z)$.
\end{theorem}

Indeed, if $M$ is homogeneous of type I and $\varphi$ is a trace the
result is immediate from Proposition~\ref{P5.4}, since the center valued
trace satisfies the conditions of $E$. The proof in the general case is
a technical adjustment to this idea.

By a similar argument, see \cite{G-N2}, one can show that if $N$ is an
injective von Neumann algebra with a normal state $\omega$, and
$(M,\varphi,\alpha)$ is a W*-dynamical system then
\begin{equation}\label{e5.6}
h_{\omega\otimes\varphi}(\id\otimes\alpha)=h_\varphi(\alpha)\;.
\end{equation}

When we apply Theorem~\ref{T5.5} to inner automorphisms we obtain

\begin{corollary}\label{C5.7}
{\rm \cite{N-S1}} \
If $(A,\varphi,\alpha)$ (resp. $(M,\varphi,\alpha)$) is a C*- (resp.
W*-) dynamical system with $A$ (resp. $M$) of type I, and $\alpha=\Ad u$
is an inner automorphism, then $h_\varphi(\alpha)=0$.
\end{corollary}

\begin{corollary}\label{C5.8}
{\rm \cite{S2}, \cite{N-S1}} \
Let $R$ be the hyperfinite II$_1$-factor, $A$ a Cartan subalgebra of
$R$ and $u$ a unitary operator in $A$. If $\varphi$ is a normal state
such that $u$ belongs to its centralizer $R_\varphi$, then
$h_\varphi(\Ad u)=0$.
\end{corollary}

\begin{proof}
As pointed out in 2.14 $A$ is conjugate to the infinite tensor product
of the diagonals $D_i$ considered in \ref{R2.19}. Thus there exists an
increasing sequence $N_1\subset N_2\subset\cdots$ of full matrix
algebras with union dense in $R$ such that $A\cong A_n\otimes B_n$,
where $A_n=N_n\cap A$, $B_n=N_n'\cap A$. Then $M_n=N_n\otimes B_n$ is of
type I and contains $A$. By Corollary~\ref{C5.7}
$h_\varphi(\Ad u|_{M_n})=0$. Since $\bigcup\limits_{n=1}^\infty M_n$ is
weakly dense in $R$,
$h_\varphi(\Ad u)=0$.
\hfill
$\Box$
\end{proof}
\bigskip

In particular it follows that if $\alpha_T$ is the automorphism of
$L^\infty(X,B,\mu)$ induced by an ergodic measure preserving
transformation $T$ and $u_T$ is the unitary in
$R=L^\infty(X,B,\mu)\times_{\alpha_T}\Z$ which implements $\alpha_T$,
and
$A$ is the masa in $R$ generated by $u_T$, then $A$ is a singular
masa, i.e. not Cartan, whenever $H(T)>0$, because by monotonicity
$H(\Ad u_T)\geq H(\alpha_T)=H(T)$.

A related consequence is due to Brown \cite{Br2}. We say a C*-algebra
$A$ is $A\T$ if it is an inductive limit of circle algebras, i.e.
$A=\overline{\bigcup\limits_{i=1}^r A_i}$, norm closure, where
$A_1\subset A_2\subset\cdots$ are C*-algebras of the form
\[
B=\bigoplus_{j=1}^k M_{n_j}(C(X_j))\;,
\]
where $X_j$ is homeomorphic to either the circle $\T$, $[0,1]$, or a
point.

\begin{corollary}\label{C5.9}
{\rm \cite{Br2}} \
If $A$ is an $A\T$ algebra and $B\subset A$ is a circle algebra we
cannot always find a sequence $(A_i)$ as above with $A_1=B$.
\end{corollary}

\begin{proof}
Let $C=\bigotimes\limits_{i\in\Z}M_i$, where $M_i=M_2(\C)$, and let
$\alpha$ be the 2-shift on $C$, and $\tau$ the unique tracial state,
see Remark~\ref{R2.19}. Then $H_\tau(\alpha)=\log 2$. Let
$A=C\times_\alpha\Z$ (see section~7 for the detailed definition), and
let $u$ be the unitary operator in $A$ which implements $\alpha$. By
\cite{BKRS} $A$ is an $A\T$ algebra, and by monotonicity $h_\tau(\Ad
u)\geq h_\tau(\alpha)=\log 2>0$. However, if $u$ belongs to a circle
algebra $A_i$ in a sequence as above, then by Corollary~\ref{C5.7}
$h_\tau(\Ad u|_{A_i})=0$, hence $h_\tau(\Ad u)=0$, which proves
that the circle algebra C*$(u)$ cannot be $A_1$ in a sequence $(A_i)$
as above.
\hfill
$\Box$
\end{proof}
\bigskip

For another entropy result on $A\T$ algebras see \cite{De}.

\section{Voiculescu's approximation entropies}

As mentioned in the introduction Voiculescu \cite{V} has introduced
entropies which are refinements of mean entropy and which provide a
very nice technique to study entropy.

Let $M$ be a hyperfinite von Neumann algebra with a faithful normal
tracial state $\tau$. Let $Pf(M)$ denote the family of finite subsets
of $M$. Modifying the notation introduced before Lemma~\ref{L2.9} we
write $\omega\subset^\delta {\mathcal X}$ if $\omega\in Pf(M)$,
${\mathcal X}\subset M$, if for each $x\in\omega$ there is a
$y\in{\mathcal X}$ such that $\|x-y\|_2<\delta$. Let further ${\mathcal
F}(M)$ denote the family of finite dimensional C*-subalgebras of $M$.
As noted in Section~2, if $A\in{\mathcal F}(M)$ then $\rank A$ is the
dimension of a masa in $A$.

\begin{definition}\label{D6.1}
{\rm \cite{V}} \
If $\omega\in Pf(M)$, $\delta>0$ put
\[
r_\tau(\omega,\delta)=\inf\{\rank
A:A\in{\mathcal F}(M), \omega\subset^\delta A\}\;,
\]
called the $\delta$-{\em rank of $\omega$}.
\end{definition}

Note that a slightly different choice for $r_\tau(\omega,\delta)$ would
be to replace $\rank A$ by $\exp(H_\tau(A))$, see \cite{C5} and
\cite{G-S2}.

\begin{definition}\label{D6.2}
{\rm \cite{V}} \
If $\alpha$ is a $\tau$-invariant automorphism of $M$ and $\delta>0$,
$\omega\in Pf(M)$ we put:
\begin{eqnarray*}
&&ha_\tau(\alpha,\omega,\delta) =
     \limsup_{n\to\infty} \frac{1}{n}
    \log r_\tau\Big( \bigcup_{j=0}^{n-1}
     \alpha^j(\omega),\delta\Big) \\
&&ha_\tau(\alpha,\omega) = \sup_{\delta>0}
          ha_\tau(\alpha,\omega,\delta) \\
&&ha_\tau(\alpha) = \sup\{ ha_\tau(\alpha,\omega):
       \omega\in Pf(M)\}\;.
\end{eqnarray*}
$ha_\tau(\alpha)$ is the {\em approximation entropy of\/} $\alpha$.
\end{definition}

An alternative is to take $\liminf$ in the definition of
$ha_\tau(\alpha,\omega,\delta)$. Then we get the {\em lower
approximation\/} entropy $\ell ha_\tau(\alpha)$.

As for the previous entropies we have
$ha_\tau(\alpha^k)=|k|ha_\tau(\alpha)$, $k\in\Z$. The proof that
$ha_\tau(\alpha^{-1})=ha_\tau(\alpha)$ is very easy; indeed
\[
r_\tau\Big( \bigcup_0^{n-1} \alpha^j(\omega),\delta\Big)
      = r_\tau(\alpha^{-n+1}\Big( \bigcup_0^{n-1}
    \alpha^j(\omega)\Big),\delta)
= r_\tau\Big( \bigcup_0^{n-1}\alpha^{-j}(\omega),\delta\Big)\;.
\]

The analogue of the Kolmogoroff-Sinai Theorem takes the following form.

\begin{proposition}\label{P6.3}
{\rm \cite{V}} \
Let $\omega_j\in Pf(M)$, $j\in\N$, $\omega_1\subset \omega_2\subset
\cdots$ be a sequence such that
$\bigcup\limits_{j\in\N}\bigcup\limits_{n\in\Z} \alpha^n(\omega_j)$
generates $M$ as a von Neumann algebra. Then
\[
ha_\tau(\alpha)=\sup_{j\in\N} ha_\tau(\alpha,\omega_j)\;.
\]
\end{proposition}

\begin{proposition}\label{P6.4}
{\rm \cite{V}} \
(i) \
If $A\in{\mathcal F}(M)$ and $\omega\in Pf(M)$ generates $A$ as a
C*-algebra then $H_\tau(A,\alpha)\leq \ell ha_\tau(\alpha,\omega)$.
\newline
(ii) \
$H(\alpha)\leq \ell ha_\tau(\alpha)\leq ha_\tau(\alpha)$.
\end{proposition}

\begin{proof}
It suffices to show (i). Let $\varepsilon>0$. By Lemma~\ref{L2.9}
there exists $\delta>0$ such that if $B\in{\mathcal F}(M)$ satisfies
$A\subset^\delta B$ then $H(A|B)<\varepsilon$. By hypothesis on
$\omega$ there exists therefore $\delta_1>0$ such that if
$\omega\subset^{\delta_1}B$ then $H(A|B)<\varepsilon$. This also
implies that if $\alpha^j(\omega)\subset^{\delta_1}B$ then
$H(\alpha^j(A)|B)<\varepsilon$. Put
$r(n)=r_\tau(\bigcup\limits_0^{n-1}\alpha^j(\omega),\delta_1)$. Then
there exists $B\in{\mathcal F}(M)$ with $\rank B=r(n)$ and
$\alpha^j(A)\subset^{\delta_1}B$ for $0\leq j\leq n-1$. Hence by
Property (F) in Section~2,
\begin{eqnarray*}
H(A,\alpha(A),\ldots,\alpha^{n-1}(A)) &\leq&
        H(B)+\sum_{j=0}^{n-1}H(\alpha^j(A)|B) \\
&\leq& \log r(n)+n\varepsilon\;,
\end{eqnarray*}
so that $H(A,\alpha)\leq ha_\tau(\alpha,\omega,\delta_1)+\varepsilon$,
proving the proposition.
\hfill
$\Box$
\end{proof}
\bigskip

In general it can be quite difficult to know when an algebra $B$ as in
the above proof satisfies $\rank B=r(n)$, hence to compute $r(n)$. A
case when it is easy is that of the $n$-shift. In the notation of
Remark~\ref{R2.19} let $R=\bigotimes\limits_{i\in\Z}(M_i,\tau_i)$ with
$M_i=M_n(\C)$, and $\alpha$ be the shift. Let
$A=M_0\in{\mathcal F}(R)$. Let, as is often done, $\omega$ be a
complete set of matrix units for $A$. By Proposition~\ref{P6.3}
$ha_\tau(\alpha,\omega)=ha_\tau(\alpha)$. On the other hand
$\bigcup\limits_{j=0}^{k-1}\alpha^j(\omega)\subset
\bigvee\limits_{j=0}^{k-1}\alpha^j(A)$, which is a I$_{n^k}$-factor,
hence
\[
r_\tau\Big( \bigcup_{j=0}^{n-1}\alpha^j(\omega),\delta\Big)
\leq n^k\qquad \mbox{for all $\delta>0$}\;.
\]
Thus, by the above and Proposition~\ref{P6.4}
$ha_\tau(\alpha)=ha_\tau(\alpha,\omega)\leq \log n=H(\alpha)\leq
ha_\tau(\alpha)$, so $ha_\tau(\alpha)=\log n$.

One test for any definition of entropy is that it should coincide
with the classical entropy in the abelian case. Via an application of
the Shannon, Breiman, McMillan Theorem the approximation entropy does
this \cite{V}.

We remarked in (\ref{2.18}) that the entropy $H(\alpha)$ is
superadditive on tensor products. For the approximation entropy the
inequality goes the other way, i.e.
\[
ha_{\tau_1\otimes \tau_2}(\alpha_1\otimes \alpha_2)\leq
ha_{\tau_1}(\alpha_1)+ha_{\tau_2}(\alpha_2)\;.
\]
Hence to show equality it suffices to show that the two entropies
coincide. In Section~10 we shall look at such cases.

In the above treatment of the approximation entropy the trace played
a minor role. If $A$ is an AF-algebra we can do essentially the same,
where we now replace the distance $\|\;\|_2$ with the operator norm.
Then we get the entropy Voiculesu denotes by hat $(\alpha)$ -- the
{\em topological approximation entropy of\/} $\alpha$.

The most flexible and therefore probably the most useful of
Voiculescu's approximation entropies are the completely positive
ones. We consider the von Neumann algebra definition first. Let
$(M,\varphi,\alpha)$ be a W*-dynamical system with $M$ injective and
$\varphi$ faithful. Let $\|x\|_\varphi=\varphi(x^\ast x)^{1/2}$ be the
$\varphi$-norm on $M$. Let
\begin{eqnarray*}
&&
CPA(M,\varphi)=\{ (\rho,\psi,B):
    B\mbox{ a finite dimensional C*-algebra, }
     \rho:M\to B,\; \psi:B\to M \\
&&\qquad\quad  \mbox{are unital completely positive maps such that }
\varphi\circ\psi\circ\rho=\varphi\}\;.
\end{eqnarray*}

\begin{definition}\label{D6.5}
{\rm \cite{V}} \
If $\omega\in Pf(M)$ and $\delta>0$ the {\em completely positive
$\delta$-rank\/} is
\begin{eqnarray*}
&&rcp_\varphi(\omega,\delta) \\
&&\qquad =\inf\{ \rank B:(\rho,\psi,B)\in CPA(M,\varphi),
    \|\psi\circ\rho(x)-x\|_\varphi<\delta\mbox{ for all
$x\in\omega$}\}\;.
\end{eqnarray*}
\end{definition}

Then we continue as in Definition~\ref{D6.2} to define the {\em
completely positive approximation entropy\/} $hcpa_\varphi(\alpha)$.
Again we can prove much the same results as for the approximation
entropy
$ha_\tau(\alpha)$.

The C*-algebra version of the above definition is like the
corresponding entropy hat $(\alpha)$ independent of invariant states.
Voiculescu defined this entropy for nuclear C*-algebras, but later on
Brown \cite{Br1} saw that one can develop the theory for exact
C*-algebras.

\begin{definition}\label{D6.6}
{\rm \cite{Br1}} \
Let $A$ be a C*-algebra and $\pi:A\to B(H)$ a faithful
$\ast$-representation. Then
\begin{eqnarray*}
&&CPA(\pi,A)=\{ (\rho,\psi,B):\rho:A\to B,\psi:B\to B(H)
      \mbox{ are contractive} \\
&&\qquad\quad \mbox{ completely positive maps, $B$ is finite
dimensional C*-algebra}\}
\end{eqnarray*}
\end{definition}

Let $\omega\in Pf(A)$, $\delta>0$. Then
\begin{eqnarray*}
&&rcp(\pi,\omega,\delta)=\inf\{ \rank B:(\rho,\psi,B)\in CPA(\pi,A),
   \mbox{ and} \\
&&\qquad \quad \|\psi\circ\rho(x)-\pi(x)\|<\delta\; \mbox{ for all
$x\in\omega$}\}\;.
\end{eqnarray*}

It follows from \cite{K} that the C*-algebras for which this definition
makes sense are the exact C*-algebras. We shall therefore assume $A$ is
exact and define the {\em topological entropy of\/} $\alpha\in\Aut(A)$,
denoted by $ht(\pi,\alpha)$ as in Definition~\ref{D6.2}.

The first result to be proved is that the definition is independent of
$\pi$, hence we can define
\[
ht(\alpha)=ht(\pi,\alpha)\;,
\]
or if $A\subset B(H)$ as $ht(\id_A,\alpha)$. The proof is a good
illustration of the techniques involved. We may assume $A$ is unital.
Let $\pi_i:A\to B(H_i)$, $i=1,2$, be faithful $\ast$-representations.
Let $\omega\in Pf(A)$, $\delta>0$. It suffices by symmetry to show
\begin{equation}\label{e6.7}
rcp(\pi_1,\omega,\delta)\geq rcp(\pi_2,\omega,\delta)\;.
\end{equation}
Choose $(\rho,\psi,B)\in CPA(\pi_1,A)$ such that $\rank
B=rcp(\pi_1,\omega,\delta)$, and
\[
\|\psi\circ\rho(x)-\pi_1(x)\|<\delta\;,\qquad x\in\omega\;.
\]
Consider the map $\pi_2\circ\pi_1^{-1}:\pi_1(A)\to B(H_2)$. From
Arveson's extension theorem for completely positive maps \cite{Ar}
there exists a unital completely positive map $\Phi:B(H_1)\to B(H_2)$
extending $\pi_2\circ\pi_1^{-1}$. Thus we have
$(\rho,\Phi\circ\psi,B)\in CPA(\pi_2,A)$ and
$\|\Phi\circ\psi\circ\rho(x)-\pi_2(x)\|<\delta$ for $x\in\omega$,
since $\pi_2(x)=\Phi\circ\pi_1(x)$. Thus (\ref{e6.7}) follows.
\hfill
$\Box$
\bigskip

Again we can prove the basic properties of entropy. Note that
monotonicity is an easy consequence of the fact that a C*-subalgebra
of an exact C*-algebra is itself exact. The analogous result is not
true for nuclear C*-algebra. We conclude this section with a theorem
which compares the entropies defined so far.

\begin{theorem}\label{T6.8}
{\rm \cite{V}} \
(i) \
If $(M,\tau,\alpha)$ is a W*-dynamical system with $\tau$ a trace then
$H_\tau(\alpha)\leq hcpa_\tau(\alpha)\leq ha_\tau(\alpha)$
\newline
(ii) \ {\rm \cite{V}} \
If $(A,\varphi,\alpha)$ is a C*-dynamical system with $A$ an
AF-algebra then $ht(\alpha)\leq hat(\alpha)$.
\newline
(iii) \
{\rm \cite{V}, \cite{D2}} \
If in (ii) $A$ is exact then $h_\varphi(\alpha)\leq ht(\alpha)$.
\end{theorem}

\section{Crossed products}

If $(A,\phi,\alpha)$ is a C*-dynamical system a natural problem is to
compute the entropy of the extension of $\alpha$ to the crossed
product $A\times_\alpha\Z$. More generally, if $G$ is a discrete
subgroup of $\Aut A$ and $\beta\in\Aut A$ commutes with $G$, compute the
entropy of the extension of $\beta$ to $A\times G$. The first positive
result is due to Voiculescu \cite{V}, who showed that for an ergodic
measure preserving Bernoulli transformation $T$ on a Lebesgue
probability space
$(X,B,\mu)$, $H(T)=H(\Ad u_T)$, where $u_T$ is the unitary operator in
$L^\infty(X,B,\mu)\times_T\Z$ which implements $T$. Later on several
extensions have appeared, see \cite{Br1}, \cite{B-C}, \cite{D-S},
\cite{G-N2}. We first recall the definition of crossed products.

Let $A$ be a unital C*-algebra, $G$ a discrete group, and
$\alpha:G\to\Aut A$ a group homomorphism. Let $\sigma:A\to B(H)$ be a
faithful nondegenerate representation. Let
\[
\pi:A\to B(\ell^2(G,H))\cong B(\ell^2(G))\otimes B(H)
\]
be the representation given by
\[
(\pi(x)\xi)(h)=\sigma(\alpha_{h^{-1}}(x))(\xi(h))\;,\qquad
x\in A\;,\quad \xi\in \ell^2(G,H)\;,\quad h\in G\;,
\]
and let $\lambda$ be the unitary representation of $G$ on
$\ell^2(G,H)$ given by
\[
(\lambda_g\xi)(h)=\xi(g^{-1}h)\;,\qquad \xi\in \ell^2(G,H)\;,\quad
g,h\in G\;.
\]
Then we have
\[
\lambda_{g^{-1}}\pi(x)\lambda_g=\pi(\alpha_g(x))\;,\qquad
x\in A\;,\quad g\in G\;.
\]
The {\em reduced crossed product C*-algebra\/} $A\times_\alpha G$ is
the norm closure of the linear span of the set
$\{\pi(x)\lambda_g:x\in A,g\in G\}$. Up to isomorphism $A\times_\alpha
G$ is independent of the choice of $\sigma$, so for simplicity we
assume henceforth that $\sigma$ is the identity map. Let
$\{\xi_h\}_{h\in G}$ be the standard orthonormal basis in $\ell^2(G)$,
so $\xi_h(g)=\delta_{g,h}$, $g,h\in G$. Then if $\xi=\xi_h\otimes\psi$
with $\psi\in H$ we have
\[
(\lambda_g\xi)(h)=\xi_{g^{-1}h}\otimes\psi=
((\ell_g\otimes 1)\xi)(h)\;,
\]
where $\ell_g$ is the left regular representation of $G$. Furthermore
\[
\pi(x)\xi=\pi(x)(\xi_h\otimes\psi)
=\xi_h\otimes\alpha_{g^{-1}}(x)\psi\;.
\]
By the above, since we may consider $A$ as a subalgebra of
$B(\ell^2(G,H))$, we may also assume from the outset that $\alpha$ is
implemented by a unitary representation $g\to U_g$, $g\in G$, of $G$.
Thus we have
\[
(1\otimes U_g)^\ast \pi(x)\lambda_h(1\otimes
U_g)=\pi(\alpha_g(x))\lambda_h\;.
\]

For simplicity let us assume $G$ is abelian -- the argument works for
$G$ amenable. We follow the approach of \cite{S-S} and \cite{Br1}. Let
$e_{p,q}\in B(\ell^2(G))$ denote the standard matrix units, i.e.
\[
e_{p,q}(\xi_t)=\delta_{q,t}\xi_p\;,
\]
where $\delta_{q,t}$ is the Kronecker $\delta$. Then we have
\[
\pi(x)\lambda_g=\sum_{t\in G} e_{t,t-g}\otimes\alpha_{-t}(x)\;,
\qquad x\in A\;,\quad g\in G\;.
\]
In particular, if $\beta\in\Aut A$, and we assume as before that
$\beta=\Ad v$ for a unitary $v\in B(H)$, then if $\beta$ commutes with
all $\alpha_g$ then
\begin{eqnarray*}
\Ad(1\otimes v)(\pi(x)\lambda_g) &=& \sum e_{t,t-g}\otimes
     v\,\alpha_{-t}(x)v^\ast = \\
&=& \sum e_{t,t-g}\otimes \alpha_{-t}(\beta(x)) \\
&=& \pi(\beta(x))\lambda_g\;.
\end{eqnarray*}
Thus $\beta$ extends to an automorphism $\hat{\beta}=\Ad(1\otimes v)$
of $A\times_\alpha G$.

If $F\subset G$ is a finite set let $P_F$ denote the orthogonal
projection of $\ell^2(G)$ onto span $\{\xi_t:t\in F\}$. Then we find
\[
(P_F\otimes 1)(\pi(x)\lambda_g)(P_F\otimes 1)=
\sum_{t\in F\cap(F+g)}e_{t,t-g}\otimes
\alpha_{-t}(x)\in M_F\otimes A\;,
\]
where $M_F=P_FB(\ell^2(G))P_F$.

In order to compute the entropy of $\hat{\beta}$ on $A\times_\alpha G$
the idea is now to start with a triple $(B,\rho,\psi)\in CPA(\id_A,A)$
and extend it to a triple $(M_F\otimes B,\Phi,\Psi)\in
CPA(\id_{A\times_\alpha G},A\times_\alpha G)$ such that we can control
the estimates. If
$f\in L^\infty(G)$ has support contained in $F$ let $m_f$ denote the
corresponding multiplication operator on $\ell^2(G)$, and define
\[
T_f(x)=\sum_{t\in G} \ell_g^\ast\otimes U_g
(m_f\otimes 1)x((m_f^\ast\otimes 1)
\ell_g\otimes U_g^\ast\;,\qquad x\in B(\ell^2(G,H))\;.
\]
Note that by amenability we can assume $\|f\|_2=1$ and
$f\ast\widetilde{f}(g_i)$ is close to 1 on a given set $g_1,\ldots,g_k$
determining $F$ where $\widetilde{f}(g)=
\overline{f(-g)}$. With $\rho$ and $\psi$ as
above we put
\[
\Phi_F(x)=(P_F\otimes 1)x P_F\otimes 1)\;,\qquad x\in
B(\ell^2(G,H))\;.
\]
Then $\Phi_F(A\times_\alpha G)\subset M_F\otimes A$, so that
\[
(M_F\otimes B,(1\otimes\rho)\circ\Phi_F,T_f\circ(1\otimes\psi))\in
CPA({\rm id}_{A\times_\alpha G},A\times_\alpha G)
\]
is the desired triple extending $(B,\rho,\psi)$. Since
$\rank(M_F\otimes B)=\card F\cdot\rank B$, all that remains is to
choose $F$ with some care depending on a given set $\omega\in
Pf(A\times_\alpha G)$, which we may suppose is of the form
$\omega=\{\pi(x_i)\lambda_{g_i};\ i=1,2,\ldots,n\}$.

The above construction essentially works for all the different
entropies considered, and even for $\beta\in\Aut A$ commuting with all
$\alpha_g$ when $G$ is amenable.

\begin{theorem}\label{T7.1}
Let $A$ be a unital C*-algebra, $G$ a discrete amenable group and
$\alpha:G\to\Aut A$ a representation. Let $\beta\in\Aut A$ commute with
all $\alpha_g$, $g\in G$. Let $\hat{\beta}$ be the natural extension of
$\beta$ to $\Aut(A\times_\alpha G)$. Then we have
\begin{itemize}
\item[(i)]
{\rm \cite{D-S}, \cite{C6}.} \
If $A$ is exact then $ht(\hat{\beta})=ht(\beta)$.

\item[(ii)]
{\rm \cite{G-N2}.} \
If $A$ is an injective von Neumann algebra and $\varphi$ a normal state
which is both $G$- and $\beta$-invariant then, if $\varphi$ is
identified with its canonical extension to $A\times_\alpha G$,
\[
hcpa_\varphi(\hat{\beta})=hcpa_\varphi(\beta)\;,
\quad \mbox{ and }\quad h_\varphi(\hat{\beta})=h_\varphi(\beta)\;.
\]
\end{itemize}
\end{theorem}

Note that when $A=L^\infty(X,B,\mu)$ and $\beta=\alpha_1=\alpha_T$,
$G=\Z$, the theorem implies the result of Voiculescu alluded to in the
first paragraph of the section. When $G$ is abelian and
$\beta=\alpha_g$ some $g\in G$, part (i) was proved by Brown
\cite{Br1}. A variation of (ii) can also be found in \cite{B-C}.

Sometimes one can prove results on operator algebras by representing
them as crossed products, see e.g. the proof of Corollary~\ref{C5.9}.
Another example is ${\mathcal O}_\infty$ -- the universal C*-algebra
generated by isometries $\{S_i\}_{i\in\Z}$ which satisfy the relation
\[
\sum_{i=-r}^r S_iS_i^\ast\leq 1\qquad \mbox{for all $\;r\in\N$}\;.
\]
Every bijection $\alpha:\Z\to\Z$ defines an automorphism, also denoted
by $\alpha$ of ${\mathcal O}_\infty$ by $\alpha(S_i)=S_{\alpha(i)}$. By
\cite{Cu} there exist an AF-algebra $B$, $\Phi\in\Aut B$, an imbedding
$\pi:{\mathcal O}_\infty\to B\times_\Phi\Z$, and a projection $p\in B$,
such that $\pi({\mathcal O}_\infty)=p(B\times_\Phi\Z)p$. By using
techniques similar to those used to prove Theorem~\ref{T7.1} we have

\begin{theorem}\label{T7.2}
{\rm \cite{B-C}} \
If $\alpha\in\Aut {\mathcal O}_\infty$ is induced by a bijective
function $\alpha:\Z\to\Z$ then $ht(\alpha)=0$. In particular, if
$\varphi$ is an $\alpha$-invariant state on ${\mathcal O}_\infty$ then
$h_\varphi(\alpha)=0$.
\end{theorem}

Note that the last statement follows from the first and
Theorem~\ref{T6.8} since $h_\varphi(\alpha)\leq ht(\alpha)$, since
${\mathcal O}_\infty$ is nuclear. For a closely related result see
\cite{C-N}. This theorem is the first we shall encounter, which shows
that if a C*-dynamical system $(A,\varphi,\alpha)$ is highly nonabelian
then the entropy of $\alpha$ tends to be small.

A problem related to the above is the computation of the entropy of the
canonical endomorphism $\Phi$ of the C*-algebra ${\mathcal O}_n$ of
Cuntz \cite{Cu}, which is the C*-algebra generated by $n$ isometries
$S_1,\ldots,S_n$ such that $\sum\limits_{i=1}^n S_iS_i^\ast=1$.
Analogously to ${\mathcal O}_\infty$, ${\mathcal O}_n$ can be written
as a crossed product $B\times_\sigma\N$, where
$B=\bigotimes\limits_{i\in\N} M_i$ with $M_i=M_n(\C)$, $\sigma$ is the
shift to the right, and $\varphi$ the canonical state extending the
trace on $B$. The {\em canonical endomorphism\/} $\Phi$ is defined by
\[
\Phi(x)=\sum_{i=1}^n S_i x S_i^\ast\;,\qquad x\in{\mathcal O}_n\;.
\]
It is a simple task to extend the entropies $h_\varphi$ and $ht$ to
endomorphisms. We have

\begin{theorem}\label{T7.3}
{\rm \cite{C4}} \
The canonical endomorphism $\Phi$ on ${\mathcal O}_n$ satisfies
\[
ht(\Phi)=h_\varphi(\Phi)=\log n\;.
\]
\end{theorem}

The result has a natural extension to the Cuntz-Krieger algebra
${\mathcal O}_A$ defined by an irreducible $n\times n$ matrix which is
not a permutation matrix. Then we have \cite{Bo-Go}
\[
ht(\alpha)=\log r(A)\;,
\]
where $r(A)$ is the spectral radius of $A$. For further extensions see
\cite{PWY}.

\section{Free products}

In Theorem~\ref{T7.2} we saw that the shift on ${\mathcal O}_\infty$ has
entropy zero. The first example of a highly nonabelian dynamical system
where the entropy is zero, was the shift on the II$_1$-factor
$L(\F_\infty)$ obtained from the left regular representation of the
free group in infinite number of generators \cite{S1}. This phenomenon
was rather surprising because the shift is so ergodic that there is no
globally invariant injective von Neumann subalgebra except for the
scalars. We shall in the present section study extensions of the above
shift.

We first recall the definitions. Let $I$ be an index set, and for each
$\iota\in I$ let $A_\iota$ be a unital C*-algebra and $\varphi_\iota$ a
state on $A_\iota$ Let $(\pi_\iota,H_\iota,\xi_\iota)$ be the
GNS-representation of $\varphi_\iota$, $\iota\in I$. Let
$H_\iota^0=H_\iota\ominus\C\xi_\iota$ and
$(H,\xi)=\mathop{\ast}\limits_{\iota\in I}(H_\iota,\xi_\iota)$ be the
free product. Put
\[
H(\iota)=\C\xi\oplus\bigoplus_{n\geq 1}
\left(
\bigoplus_{\iota_1\not=\iota_2\not=\cdots\not=\iota_n\atop
\iota_1\not= \iota}H_{\iota_1}^0\otimes\cdots\otimes H_{\iota_n}^0
\right)
\]
We have unitary operators $V_\iota:H_\iota\otimes H(\iota)\to H$
defined by
\begin{eqnarray*}
&&\xi_\iota\otimes\xi\to\xi \\
&&H_\iota^0\otimes\xi\to H_\iota^0\qquad \mbox{by}\quad
    \eta\otimes\xi\to\eta \\
&&\xi_\iota\otimes(H_{\iota_1}^0\otimes\cdots\otimes H_{\iota_n}^0)
    \to H_{\iota_1}^0\otimes\cdots\otimes H_{\iota_n}^0\qquad
    \mbox{by}\quad \xi_\iota\otimes\eta\to\eta,\;\,
    \iota_1\not=\iota \\
&&H_\iota^0\otimes(H_{\iota_1}^0\otimes\cdots\otimes H_{\iota_n}^0)
    \to H_\iota^0\otimes H_{\iota_1}^0\otimes\cdots\otimes
H_{\iota_n}^0\qquad
    \mbox{by}\quad \psi\otimes\eta\to\psi\otimes\eta,\;\,
    \iota_1\not=\iota\;.
\end{eqnarray*}
The representation $\lambda_\iota:A_\iota\to B(H)$ is defined by
\[
\lambda_\iota(x)=V_\iota(\pi_\iota(x)\otimes
1_{H(\iota)})V_\iota^\ast\;,\qquad x\in A_\iota\;.
\]
The free product representation
$\pi=\mathop{\ast}\limits_{\iota\in I}\pi_\iota: \ast A_\iota\to
B(H)$ is the *-homomorphism of the free product C*-algebra $(\ast
A_\iota,\ast \lambda_\iota)\to B(H)$, using the universal property of
the free product. When we write $(A,\varphi)=(\ast
A_\iota,\ast\varphi_\iota)_{\iota\in I}$ we shall mean $\ast A_\iota$
in the representation $\pi$, i.e. we shall mean $\pi(\ast
A_\iota)\subset B(H)$.

There are now two approaches; the first is to show that $ht(\alpha)$
for $\alpha$ the shift, or rather $\alpha$ defined via a bijection
$\Z\to\Z$ like for ${\mathcal O}_\infty$ in Theorem~\ref{T7.2}, and use
the inequalities in Theorem~\ref{T6.8} to conclude that the other
entropies are zero. This works if the $A_\iota$ are all exact, because
then $A$ is exact \cite{D1}.

\begin{theorem} \label{T8.1}
{\rm \cite{D2}} \
In the above notation assume each $A_\iota$ is finite dimensional and
let $\sigma$ be a permutation of $I$ such that for all $\iota\in I$
there is an isomorphism $\alpha_\iota:A_\iota\to A_{\sigma(\iota)}$
such that $\varphi_{\sigma(\iota)}\circ\alpha=\varphi_\iota$. Then
there exists a unique $\alpha\in\Aut A$ such that
$\alpha\circ\lambda_\iota=\lambda_{\sigma(\iota)}\circ\alpha_\iota$,
$\iota\in I$, and $ht(\alpha)=0$.
\end{theorem}

Dykema proves a more general result than the above. The theorem is not
as restricted as it looks, because if we let $B_\iota=A_{2\iota}\ast
A_{2\iota+1}$ we can write $A$ as a free product $\ast B_\iota$ of a
large class of C*-algebras. For example in this way we can imbed
${\mathcal O}_\infty$ in a free product to recover Theorem~\ref{T7.2}.
Similarly we can obtain the announced result on the shift of
$C_r^\ast(\F_\infty)$, and hence by (\ref{e3.8}) of
$L^\infty(\F_\infty)$, from the inequality $h_\varphi(\alpha)\leq
ht(\alpha)$.

For the other approach we note that we have no Kolmogoroff-Sinai
Theorem in the nonnuclear case, so we must go directly at the
definition of $h_\varphi(\gamma_1,\ldots,\gamma_n)$. Recall that we
then considered a positive map $P:A\to B$, where $B$ is a finite
dimensional abelian C*-algebra with a state $\mu$ such that $\mu\circ
P=\varphi$. Denote by $\|x\|_\mu=\mu(x^\ast x)^{1/2}$.

\begin{lemma}\label{L8.2}
{\rm \cite{S3}} \
Let $(A,\varphi)=(\ast A_\iota,\ast\varphi_\iota)_{\iota\in I}$ be a
free product of unital C*-algebras. Let $B$ be an abelian C*-algebra
with a state $\mu$. Suppose $P:A\to B$ is a unital positive linear map
such that $\mu\circ P=\varphi$. Then given $\varepsilon>0$ there is
$J\subset I$ with $\card J\leq\left[\frac{1000}{\varepsilon}\right]+1$
such that
\[
\|P(x)-\varphi(x)1\|_\mu<\varepsilon \|x\|\;,\qquad
x\in A_\iota\;,\quad \iota\not\in J\;.
\]
\end{lemma}

Thus $P$ is essentially almost constant outside the subalgebra
$\mathop{\ast}\limits_{\iota\in J}A_\iota$. If we now assume
$A_\iota=A_0$, $\iota\in I=\Z$ and letting for example $\alpha$ be the
free shift on $A$ arising from the shift $\iota\to\iota+1$ on $\Z$,
then it is not hard to go through the different steps in the definition
of the CNT-entropy $h_\varphi(\alpha)$ to conclude:

\begin{theorem}\label{T8.3}
{\rm \cite{S3}} \
If $A_0$ is a unital C*-algebra and $\varphi_0$ a state on $A_0$, and
$A_i=A_0$, $\varphi_i=\varphi_0$, $i\in\Z$, then the free shift
$\alpha$ on $(A,\varphi)=(\ast A_i,\ast\varphi_i)_{i\in\Z}$ has entropy
$h_\varphi(\alpha)=0$.
\end{theorem}

\begin{remark}\label{R8.4}
\rm \
It should be noted that by a result of Avitzour \cite{Av} it follows
that every stationary coupling extending $\varphi$ is of the form
$\lambda=\varphi\otimes\mu$, hence the entropy $h_\varphi'(\alpha)$
of Sauvageot and Thouvenot is zero. One can further show \cite{C3}
that if $(C,\rho,\gamma)$ is a C*-dynamical system then
\[
h_{\varphi\ast\rho}'(\alpha\ast\gamma)=h_\rho'(\gamma)\;.
\]
Similar results hold for other entropies, see \cite{C6}.
\end{remark}

\section{Binary shifts}

A rich class of C*-dynamical systems is obtained from bitstreams, i.e.
sequences $(x_n)_{n\in\N}$ with $x_n\in\{0,1\}$. Denote by $X$ the
subset of $\N$, $X=\{n\in\N: x_n=1\}$. We can construct a sequence
$(s_n)_{n\in\N}$ of symmetries, i.e. self-adjoint unitary operators on a
Hilbert space which satisfy the commutation relations
\[
s_is_j= \begin{cases}
s_js_i & \mbox{if}\quad |i-j|\not\in X,\quad
       \mbox{i.e.} \;x_{|i-j|}=0 \\
-s_js_i & \mbox{if}\quad |i-j|\in X,\quad \mbox{i.e.} \;x_{|i-j|}=1,
\end{cases}
\]
see e.g. \cite{Vi}.

Let $A(X)$ denote the C*-algebra generated by the $s_n$, $n\in\Z$. The
canonical trace $\tau$ on $A(X)$ is the one which takes the value zero
on all products $s_{i_1}\ldots s_{i_n}$, where $i_1<i_2<\cdots<i_n$,
and $\tau(1)=1$. We denote by $\alpha$ the shift automorphism of $A(X)$
defined by $\alpha(s_i)=s_{i+1}$. Then $(A(X),\tau,\alpha)$ is a
C*-dynamical system. Well-known situations from both C*-algebras and
the classical case are represented as special cases, e.g.
asymptotically abelian, proximally asymptotically abelian, K-systems,
and completely positive entropy, see \cite{G-S1}. We shall assume we
are in the nontrivial case when then set $-X\cup\{0\}\cup X$ is
nonperiodic. Let $A_n=C^\ast(s_0,\ldots,s_{n-1})$ be the C*-algebra
generated by $s_0,\ldots,s_{n-1}$, so that
\[
A_n=\bigvee_0^{n-1} \alpha^i(C^\ast(s_0))\;.
\]
We list some properties of $A_n$ and $A(X)$ which will be used in the
sequel, see \cite{E}, \cite{Po-Pr} or \cite{Vi}. Denote by $Z_n$ the
center of $A_n$. Then we have:
\begin{equation}\label{e9.1}
\mbox{There are $c_n,d_n\in\N\cup\{0\}$ such that
     $n=2d_n+c_n$}\,,\quad
A_n\cong M_{2^{d_n}}(\C)\otimes Z_n
\end{equation}
where $Z_n\cong C(\{0,1\}^{c_n})$.
\begin{eqnarray}
&&\mbox{If $e$ is a minimal projection in $Z_n$ then
    $\tau(e)=2^{-c_n}$}.
   \label{e9.2} \\
&&\mbox{There is a sequence $(m_i)$ in $\N$ such that $(c_n)$ consists
    of the concatenation}
    \label{e9.3} \\
&&\mbox{of the strings
    $(1,\ldots,m_i-1,m_i,m_i-1,\ldots,1,0)$.} \nonumber
\end{eqnarray}
In particular $c_n=0$ for an infinite number of $n$'s, hence
$A_n=M_{2^{n/2}}(\C)$ for these $n$'s, and so $A(X)$ is the CAR-algebra.
\begin{equation}\label{e9.4}
H(A_n)=\log\rank A_n=(c_n+d_n)\log 2\;.
\end{equation}
By (\ref{e9.1}) $2d_n\leq n\leq 2d_n+2c_n$, hence
\[
\liminf_n\frac{1}{n}H(A_n)=\liminf_n\frac{1}{n}
(c_n+d_n)\log 2\leq{\textstyle\frac{1}{2}}\log 2\;.
\]
Since also $\frac{1}{n}(c_n+d_n)\geq\frac{1}{2}$ we find
\begin{eqnarray}
&&\liminf_n\frac{1}{n}H(A_n)={\textstyle\frac{1}{2}}\log 2
     \label{e9.5} \\
&&\lim_{n\to\infty}\frac{c_n}{n}=0\quad
   \mbox{if and only if $\;\lim\frac{1}{n}H(A_n)
    =\frac{1}{2}\log 2$}\;.
   \label{e9.6}
\end{eqnarray}
Indeed, if $\frac{c_n}{n}\to0$ then $\frac{d_n}{n}\to\frac{1}{2}$,
hence by (\ref{e9.4}) $\frac{1}{n}H(A_n)\to\frac{1}{2}\log 2$.
Conversely, if $\lim\limits_n\frac{1}{n}H(A_n)=\frac{1}{2}\log 2$ then
by (\ref{e9.4}) $\frac{1}{n}(c_n+d_n)\to\frac{1}{2}$, hence by
(\ref{e9.1}) $\frac{c_n}{n}\to 0$.

Since $A(X)$ is the CAR-algebra the weak closure of its image in the
GNS-representation of $\tau$ is the hyperfinite II$_1$-factor $R$. When
we in the sequel consider the approximation entropies of $\alpha$
defined in Definition~\ref{D6.2} it is really the extension of $\alpha$
to $R$ that we consider. While $ha_\tau$ is subadditive on tensor
products the lower approximation entropy $\ell ha_\tau$ only satisfies
\[
\ell ha_{\tau\otimes\tau}(\alpha\otimes\alpha)\leq2\ell
ha_\tau(\alpha)\;.
\]

\begin{lemma}\label{L9.7}
{\rm \cite{NST}, \cite{G-S2}} \
With $\alpha$ and $X$ as before we have:
\begin{itemize}
\item[(i)]
$h_{\tau\otimes\tau}(\alpha\otimes\alpha)=\log 2$

\item[(ii)]
$\ell ha_\tau(\alpha)=\frac{1}{2}\log 2$
\end{itemize}
\end{lemma}

\begin{proof}
Let $A_0$ denote the C*-subalgebra of $A(X)\otimes A(X)$ generated by
the symmetries $s_i\otimes s_i$, $i\in\Z$. Then $A_0$ is abelian, and
$\tau\otimes\tau$ vanishes on each $s_i\otimes s_i$. Thus the
C*-dynamical system
$(A_0,\tau\otimes\tau,\alpha\otimes\alpha)$ is isomorphic to the
2-shift, hence has entropy $\log 2$, hence by monotonicity
$h_{\tau\otimes\tau}(\alpha\otimes\alpha)\geq\log 2$. Thus by
Proposition~\ref{P6.4}, and the inequality preceeding the lemma
\[
2\ell ha_\tau(\alpha)\geq \ell
ha_{\tau\otimes\tau}(\alpha\otimes\alpha)\geq
h_{\tau\otimes\tau}(\alpha\otimes\alpha)\geq\log 2\;.
\]
The converse inequality follows from easy estimates using (\ref{e9.5}).
\hfill
$\Box$
\end{proof}
\bigskip

>From the above lemma $h_\tau(\alpha)\in[0,\frac{1}{2}\log 2]$. In many
cases $h_\tau(\alpha)=\frac{1}{2}\log 2$.

\begin{theorem}\label{T9.8}
{\rm \cite{C1}, \cite{G-S1}} \
Suppose $X$ satisfies one of the following: $X$ is finite,
$\N\setminus X$ is finite, $X$ is contained in the even or odd
numbers, $X$ contains the odd numbers. Then
$h_\tau(\alpha)=\frac{1}{2}\log 2$.
\end{theorem}

Consider the case when $X\supset\{1,3,5,\ldots\}$. Let
$t_j=s_{2j-1}s_{2j}$, $j\in\Z$. Then the $t_j$ all commute, and as in
the proof of Lemma~\ref{L9.7} $\alpha^2$ acts as a 2-shift on the
C*-algebra they generate. Thus
\[
h_\tau(\alpha)={\textstyle\frac{1}{2}}h_\tau(\alpha^2)\geq
{\textstyle\frac{1}{2}}\log 2\;,
\]
and the conclusion follows.
\hfill
$\Box$
\bigskip

For other examples when $h_\tau(\alpha)=\frac{1}{2}\log 2$ see
\cite{Pr}. If $J\subset\N$ is finite, $J=\{i_1,\ldots,i_n:
i_1<i_1<\cdots<i_n\}$ then the operator $s_{i_1}\ldots s_{i_n}$ or
$is_{i_1}\ldots s_{i_n}$ is self-adjoint. Denote the self-adjoint
operator by $s_J$. Then $\{\alpha^j(s_J)\}_{j\in \Z}$ is a sequence of
symmetries which either commute or anticommute. Let
\[
X(J)=\{j\in\N: \alpha^j(s_J)s_J+s_J \alpha^j(s_J)=0\}
\]
be the corresponding subset of $\N$. If $I\subset\N$ we denote by
$I-I=\{n-m:m,n\in I\}$.

\begin{theorem}\label{T9.9}
{\rm \cite{NST}} \
(i) \ Assume for each finite subset $J\subset\N$ there exists an
infinite subset $I\subset\N$ such that $(I-I)\cap\N\subset X(J)$. Then
$h_\tau(\alpha)=0$.

(ii) \
There exists $X\subset\N$ such that (i) holds.
\end{theorem}

The proof of (i) consists of showing that the Sauvageot-Thouvenot
entropy $h'_\tau(\alpha)=0$. This is done by showing that a stationary
coupling $\lambda$ corresponding to an abelian system $(B,\mu,\beta)$
necessarily is of the form $\lambda=\tau\otimes\mu$.

If we combine Theorem~\ref{T9.9} with Lemma~\ref{L9.7} we have

\begin{corollary}\label{C9.10}
{\rm \cite{NST}} \
There exists $X\subset\N$ such that
$h_{\tau\otimes\tau}(\alpha\otimes\alpha)=\log 2$ while
$h_\tau(\alpha)=0$.
\end{corollary}

We leave it as open
problems whether there exists $X\subset\N$ such that
$0<h_\tau(\alpha)<\frac{1}{2}\log 2$, and whether there exists
$X\subset\N$ such that, cf. Lemma~\ref{L9.7} (ii),
$ha_\tau(\alpha)>\frac{1}{2}\log 2$.

Another example of a C*-dynamical system $(A,\tau,\alpha)$ for which
$h_{\tau\otimes\tau}(\alpha\otimes\alpha)>
h_\tau(\alpha)+h_\tau(\alpha)=0$
has been exhibited by Sauvageot \cite{Sa}, see also \cite{N-T2}.

Let $\theta\in\R$ and $A_\theta$ be the C*-algebra
generated by two unitaries $U$ and $V$ such that
\[
VU=e^{2\pi i\theta}UV\;,
\]
so $A_\theta$ is the irrational rotation algebra when $\theta$ is
irrational. If $\mu=(m,n)\in\Z^2$ let
\[
W_\theta(\mu)=e^{i\pi\theta mn}U^mV^n\;.
\]
Then linear combinations of the $W_\theta(\mu)$, $\mu\in\Z^2$,
are dense in $A_\theta$. In analogy with Bogoliubov automorphisms
each matrix $A\in SL(2,\Z)$ defines an automorphism $\sigma_A$ of
$A_\theta$ by
\[
\sigma_A(W_\theta(\mu))=W_\theta(A\mu)\;.
\]
$A_\theta$ has a canonical trace $\tau_\theta$ such that
$\tau_\theta(W_\theta(\mu))=0$ whenever $\mu\not=(0,0)$. We
make the assumption that $A$ has two real eigenvalues $\lambda$ and
$\lambda^{-1}$ with $|\lambda|>1$. Then we have

\begin{theorem}\label{T9.5}
{\rm \cite{Sa}, \cite{N-T2}} \
With the above assumption there is a subset $\Omega\subset\R$ for which
the complement $\Omega^c$ has Lebesgue measure zero and
$\Q+\lambda\Q\subset\Omega^c$, with the following properties:
\begin{itemize}
\item[(i)]
If $\theta\in\Q+\lambda\Q$ then $h_{\tau_\theta}(\sigma_A)>0\,$.

\item[(ii)]
If $\theta\in\Omega$ then $\tau_\theta$ is the unique
invariant state and
$h_{\tau_\theta}(\sigma_A)=0=h_{\tau_{-\theta}}(\sigma_A)\,$.
\end{itemize}

\noindent
Furthermore, in case (ii)
$h_{\tau_\theta\otimes\tau_{-\theta}}
(\sigma_A\otimes\sigma_A)>0$.
\end{theorem}

\section{Generators}

Corollary~\ref{C9.10} shows that the tensor product formula
$h_{\varphi\otimes\varphi}(\alpha\otimes\beta)=
h_\varphi(\alpha)+h_\varphi(\beta)$ can only hold in special cases. By
the superadditivity of $h_\varphi$, see (\ref{e3.9}), and the
subadditivity of the approximation entropy $ha_\tau$, the tensor
product formula holds whenever the two entropies coincide. Since the
latter is a refined version of mean entropy we may therefore expect the
tensor product formula to hold whenever the CNT-entropy $h_\tau$ is a
mean entropy.

In the classical case of a probability space $(X,{\mathcal B},\mu)$ with
a measure preserving nonsingular transformation $T$, a partition $P$ is
called a generator if the $\sigma$-algebra
$\bigvee\limits_{i\in\Z}T^{-i}P={\mathcal B}$, or equivalently, if $A$
is the C*-algebra generated by the atoms in $P$, then
$\bigvee\limits_{i\in\Z}\alpha^i(A)=L^\infty(X,{\mathcal B},\mu)$. In
\cite{G-S2} different candidates for nonabelian generators were
considered.

\begin{definition}\label{D10.1}
Let $M$ be a hyperfinite von Neumann algebra with a faithful normal
tracial state $\tau$. Let $\alpha$ be a $\tau$-invariant automorphism.
A finite dimensional von Neumann subalgebra $N$ of $M$ is a {\em
generator\/} for $\alpha$ if
\medskip

(i) \
$\bigvee\limits_{i\in\Z}\alpha^i(N)=M$.

(ii) \
$\bigvee\limits_{i=m}^n \alpha^i(N)$ is finite dimensional whenever
$m<n$, $m,n\in\Z$.

(iii) \
$H(N,\alpha)=\limsup\limits_n\frac{1}{n}H\left(
\bigvee\limits_0^{n-1}\alpha^i(N)\right)$.
\medskip

\noindent
If (iii) is replaced by

(iv) \
$H(\alpha)=\limsup\limits_n\frac{1}{n}H
\left( \bigvee\limits_0^{n-1}\alpha^i(N)\right)$.
\medskip

\noindent
then $N$ is called a {\em mean generator}.
\end{definition}

If $N$ is a generator then $N$ is a mean generator, and
\[
H(\alpha)=H(N,\alpha)=\lim_n\frac{1}{n}H\left(
\bigvee_0^{n-1}\alpha^i(N)\right).
\]

We say $N$ is a {\em lower generator\/} (resp. {\em lower mean
generator\/}) if we replace $\limsup\limits_n$ by $\liminf\limits_n$ in
(iii) (resp. (iv)). Recall from Definition~\ref{D6.2} that if we replace
$\rank A$ in the definition of $ha_\tau(\alpha)$ and $\ell
ha_\tau(\alpha)$ by $\exp(H(A))$ we obtained two entropies we denote
by $Ha_\tau(\alpha)$ and $\ell Ha_\tau(\alpha)$. An easy consequence of
the definitions is then

\begin{proposition}\label{P10.2}
Let $N\subset M$ be finite dimensional. If $N$ is a mean generator
(resp. lower mean generator) then $H(\alpha)=Ha_\tau(\alpha)$ (resp.
$H(\alpha)=\ell Ha_\tau(\alpha)$).
\end{proposition}

\begin{corollary}\label{C10.3}
Let $(M_i,\tau_i,\alpha_i)$ be W*-dynamical systems as above, $i=1,2$.
If
$\alpha_1$ and $\alpha_2$ have mean generators then
\[
H_{\tau_1\otimes\tau_2}(\alpha_1\otimes\alpha_2)=
H_{\tau_1}(\alpha_1)+H_{\tau_2}(\alpha_2)\;.
\]
\end{corollary}

Since the CNT-entropy remains the same when we imbed C*-dynamical
system into the W*-dynamical system obtained via the GNS-representation
due to the invariant state (\ref{e3.8}) our definitions clearly make
sense for AF-algebras.

Several well-known examples of C*-dynamical systems have generators. We
list a few.

\setcounter{subsection}{3}

\subsection{Temperley-Lieb algebras}

Let $(e_i)_{i\in\Z}$ be a sequence of projections with the porperties.
\begin{itemize}
\item[(a)]
$e_ie_{i\pm1}e_i=\lambda e_i$ for some $\lambda\in(0,\frac{1}{4}]\cup
\{\frac{1}{4}{\rm sec}^2(\frac{\pi}{m}):m\geq3\}$

\item[(b)]
$e_ie_j=e_je_i$ if $|i-j|\geq2$.

\item[(c)]
$\tau(\omega e_i)=\lambda\tau(\omega)$ if $\omega$ is a word in 1 and
$e_j$, $j<i$.
\end{itemize}

\noindent
Then the von Neumann algebra $R$ generated by the $e_i$'s is the
hyperfinite II$_1$-factor. The shift automorphism $\alpha_\lambda$
determined by $\alpha_\lambda(e_i)=e_{i+1}$ has $C^\ast(e_0)$ as a mean
generator, see \cite{G-S2}.

\subsection{Noncommutative Bernoulli shifts}

In the notation of (2.14) assume $d=2$. Let $N$ denote the centralizer
of $\varphi$ in $M_0\otimes M_1$. Then $N$ is a generator for $\alpha$
\cite{G-S2}.

\subsection{Binary shifts}

Assume as before $-X\cup\{0\}\cup X$ is nonperiodic. By
Theorem~\ref{T9.9}
and
Corollary~\ref{C10.3} we do not in general have generators for binary
shifts. In the notation of (\ref{e9.1})--(\ref{e9.6}) the following
hold.
\begin{itemize}
\item[(i)]
If $H(\alpha)=\frac{1}{2}\log 2$ then $A_1$ is a lower mean generator.
If moreover $\lim\limits_n\frac{c_n}{n}=0$ then $A_1$ is a mean
generator.

\item[(ii)]
If $X$ is contained in the even numbers then $A_2$ is a lower
generator. If moreover $\lim\limits_n\frac{c_n}{n}=0$ then $A_2$ is a
generator.

\item[(iii)]
If $X$ is contained in the odd numbers then $A_1$ is a lower generator.
If moreover $\lim\limits_n\frac{c_n}{n}=0$ then $A_1$ is a generator.
\end{itemize}

There are many ways to generalize the concept of generators. Instead of
looking at $\bigvee\limits_0^{n-1}\alpha^i(N)$ we can consider an
increasing sequence of finite dimensional algebras and endomorphisms
instead of automorphisms.

\setcounter{theorem}{6}

\begin{definition}\label{D10.7}
{\rm \cite{S4}} \
Let $M$ be a hyperfinite von Neumann algebra with a faithful normal
tracial state. Let $\alpha$ be a $\tau$-invariant endomorphism. An
increasing sequence $(A_n)_{n\in\N}$ of finite dimensional von Neumann
subalgebras of $M$ is a {\em generating sequence\/} for $\alpha$ if

(i) \ $\bigvee\limits_{n=1}^\infty A_n=M$,

(ii) \ $\alpha(A_n)\subset A_{n+1}$, $n\in\N$,

(iii) \ $H(\alpha)=\lim\limits_n\frac{1}{n}H(A_n)$.
\end{definition}

By (iii) it is clear that the tensor product formula holds in the
presence of generating sequences.

We say a sequence $(A_n)_{n\in\N}$ such that (i) and (ii) hold,
satisfies the {\em commuting square condition\/} if the diagram
\[
\xymatrix@1@=5pt{A_{n+1} & {\subset} & M \\
{\bigcup} & & {\bigcup} \\
\alpha(A_n) & \subset & \alpha(A_{n+1})
}
\]
is a commuting square for all $n\in\N$, i.e.
$E_{\alpha(A_n)}=E_{\alpha(A_{n+1})}\circ E_{A_{n+1}}$.

With the above assumptions we can generalize the classical result that
if $P$ is a generator for a measure preserving nonsingular
transformation $T$ on  probability space $(X,{\mathcal B},\mu)$ then
$H(T)$ is given by relative entropy,
\[
H(T)=\lim_{n\to\infty} H\left( \bigvee_0^n T^{-i}(P)\Big|
\bigvee_1^n T^{-i}(P)\right).
\]
In our case relative entropy is defined in Definition~\ref{D2.8}. If
$N$ is a von Neumann algebra we denote by $Z(N)$ the center of $N$.

\begin{theorem}\label{T10.5}
{\rm \cite{S4}} \
Let $M,\tau,\alpha$ be as above and suppose $H(\alpha)<\omega$. Suppose
$(A_n)_{n\in\N\cup\{0\}}$ is a generating sequence for $\alpha$
satisfying the commuting square condition. Then we have:

(i) \
$\lim\limits_{n\to\infty}\frac{1}{n}H(Z(A_n))$ exists.

(ii) \
$H(\alpha)=\frac{1}{2}H(M| \alpha(M))+
\frac{1}{2}\lim\limits_{n\to\infty}\frac{1}{n}H(Z(A_n))$.

\noindent
Furthermore, if $M$ is of type I then $H(\alpha)=H(M|\alpha(M))$.
\end{theorem}

Pimsner and Popa \cite{P-P} found an explicit formula for the relative
entropy $H(P|Q)$ when $P\supset Q$ are finite dimensional C*-algebras.
For the sequence $(A_n)$ the formula becomes
\begin{eqnarray}\label{e10.9}
&&H(A_n)|\alpha(A_{n-1})) \\
&&\qquad =
    2(H(A_n)-H(\alpha(A_{n-1}))
    -(H(Z(A_n))-H(Z(\alpha(A_{n-1})))+c_n\;, \nonumber
\end{eqnarray}
where $c_n$ is a real number depending on the multiplicity of the
imbedding $\alpha(A_{n-1})\subset A_n$. Since
$H(\alpha(A_{n-1}))=H(A_{n-1})$ and similarly for the centers, if we
add the equations in (\ref{e10.9}) we find,
\begin{eqnarray*}
&&\frac{1}{k}\sum_{n=1}^k H(A_n| \alpha(A_{n-1})) \\
&&\qquad\qquad
    =\frac{2}{k}H(A_k)-\frac{2}{k}H(A_0)
-\frac{1}{k}H(Z(A_k))
+\frac{1}{k}H(Z(A_0))+\frac{1}{k}\sum_{n=1}^k c_n\;.
\end{eqnarray*}
Since $(A_n)$ satisfies the commuting square condition the left side
converges to $H(M|\alpha(M))$, \cite{P-P}. Some analysis shows
$\frac{1}{k}\sum\limits_i^k c_n\to 0$. Since the terms involving $A_0$
go to zero and $\frac{2}{k}H(A_k)\to 2H(\alpha)$, the proof is
completed.
\hfill
$\Box$
\bigskip

If $R$ is the hyperfinite II$_1$-factor then by a formula of Pimsner
and Popa \cite{P-P} we get

\begin{corollary}\label{C10.10}
{\rm \cite{S4}} \
If in Theorem~\ref{T10.5}
$M=R$ is the hyperfinite II$_1$-factor we have

(i) \ $\lim\limits_n \frac{1}{n}H(Z(A_n))$ exists.

(ii) \
$R\cap \alpha(R)'$ is atomic with minimal projections $f_k$,
$\sum\limits_k f_k=1$,

(iii) \
$H(\alpha)=H(R\cap \alpha(R)')
+\frac{1}{2}\sum\limits_k \tau(f_k)\log[R_{f_k}: \alpha(R)_{f_k}]
+\frac{1}{2}\lim\limits_n \frac{1}{n}H(Z(A_n))$.
\end{corollary}

Here $R_f$ means the $R$ cut down by $f$, and $[P:Q]$ is the Jones
index \cite{J}. In particular, if $R\cap \alpha(R)'=\C$ then
\begin{equation}\label{e10.11}
H(\alpha)={\textstyle\frac{1}{2}}\log[R:\alpha(R)]
+{\textstyle\frac{1}{2}}\lim_n \frac{1}{n}H(Z(A_n))\;.
\end{equation}

The natural sequences considered in Examples (10.4)--(10.6) all satisfy
the commuting square conditions, and in (10.4) and (10.6)
(\ref{e10.11}) reduces to
\begin{equation}\label{e10.12}
H(\alpha)={\textstyle\frac{1}{2}}\log[R:\alpha(R)]\;,
\end{equation}
a formula which is well-known in those cases, see \cite{J},
\cite{P-P}, \cite{P}.

If $N\subset M$ is an inclusion of type III factors with a common
separating and cyclic vector $\xi$ the conjugations $J_M$ and $J_N$
defined by $\xi$ define an endomorphism $\gamma$ of $M$ into $N$ by
\[
\gamma(x)=J_NJ_M x J_MJ_N\;,\qquad x\in M\;,
\]
see \cite{L}. If the index is finite Choda \cite{C2} showed a formula
like (\ref{e10.12}) for the entropy of $\gamma$ with respect to a
natural invariant state.

The II$_1$-analogue $\Gamma$ of $\gamma$ is defined for an irreducible
inclusion of II$_1$-factors $N\subset M$ of finite index. We then get a
tower
\[
N=M_{-1}\subset M=M_0\subset M_1\subset \cdots
\]
The trace $\tau$ of $M$ is extended naturally to a trace, also denoted
by $\tau$, on each $M_k$. The canonical conjugation $J_k$ on the
Hilbert space $L^2(M_k,\tau)$ is defined by $J_kx=x^\ast$, $x\in M_k$.
We denote by $M_\infty$ the closure of $\bigcup\limits_k M_k$ in the
GNS-representation of $\tau$ and consider $M$ as a subfactor of
$M_\infty$. Let $R=M'\cap M_\infty$. Then $R$ is the hyperfinite
II$_1$-factor, and the canonical shift $\Gamma$ on $R$ is defined by
\[
\Gamma(x)=J_{n+1}J_n x J_nJ_{n+1}
\]
for $x\in M'\cap M_{2k}$, $n\geq k$. This definition does not depend on
$n$. Let $A_n=M'\cap M_{2n}$. Then the sequence $(A_n)$ is a generating
sequence which satisfies the commuting square condition. It follows
from Corollary~10.9 that $H(\Gamma)$ is given by
(\ref{e10.11}). This formula was first shown by Hiai \cite{H}. With
certain extra assumptions he and Choda \cite{C2}, \cite{C-H} showed
that $H(\Gamma)$ is given by (\ref{e10.12}).

\section{The variational principle}

The variational principle appears as an important ingredient both in
classical ergodic theory and in spin lattice systems in the C*-algebra
formalism of quantum statistical mechanics. The principle can be well
described in the finite dimensional case. Indeed, let $B$ be a finite
dimensional C*-algebra and $\Tr_B$ the canonical trace on $B$, i.e.
$\Tr_B(e)=1$ for all minimal projections $e\in B$. If $\varphi$ is a
state on $B$ let $K_\varphi\in B^+$ be its density operator, so
$\varphi(x)=\Tr_B(K_\varphi x)$. Then the mean entropy of $\varphi$ is
\[
S(\varphi)=\Tr_B(\eta(K_\varphi))\;.
\]
The variational principle takes the form: If $H\in B_{\sa}$ then
\begin{equation}\label{e11.1}
S(\varphi)-\varphi(H)\leq \log \Tr_B(e^{-H})\;,
\end{equation}
with equality if and only if
\[
K_\varphi=\frac{e^{-H}}{\Tr_B(e^{-H})}\;.
\]
In this case $\varphi$ is called the Gibbs state. Note that $\varphi$
is then a KMS-state at $\beta=1$ for the one-parameter group
\[
\sigma_t^H(x)=e^{-\ita H} x e^{\ita H}\;,
\]
where we recall that a state $\varphi$ on a C*-algebra $A$ is a
KMS-state for a one-parameter group $\sigma_t$ at temperature $\beta$ if
\[
\varphi(ab)=\varphi(b\sigma_{i\beta}(a))
\]
for all analytic elements $a,b\in A$, see \cite[Ch.~5]{B-R}.

The above variational principle has been extended to automorphisms of
$C(X)$ when $X$ is a compact metric space, see \cite[Ch.~9]{W}, and to
spin lattice systems, see \cite[Ch.~6]{B-R}. In the latter case the
C*-algebra is an infinite tensor product indexed by $\Z^\nu$,
$\nu\in\N$, and the automorphism group consists of the shifts. In the
years around 1970 the variational principle was intensely studied by
mathematical physicists, see the comments to Section~6.2.4 in
\cite{B-R}, and solved in a way naturally extending the results in the
finite dimensional case. The entropy used, was mean entropy. However,
Moriya \cite{M} has shown that one gets the same results with
CNT-entropy replacing mean entropy.

The C*-dynamical systems $(A,\alpha)$ considered above are all
asymptotically abelian, i.e. $\|[\alpha^n(x),y]\|\to0$ as $n\to\infty$
for all $x,y\in A$. Following work of Neshveyev and the author
\cite{N-S2} we shall in the present section see that the variational
principle has a natural extension to a class of asymptotically abelian
C*-algebras. For another approach in the case of Cuntz-Krieger algebras
see \cite{PWY}.
First we shall need to extend the right side of
(\ref{e11.1}) to C*-algebras. We modify the definition of topological
entropy in Definition~\ref{D6.6}.

Let $A$ be a nuclear C*-algebra with unit and $\alpha$ an automorphism.
Let
\begin{eqnarray*}
&&{\rm CPA}(A)=\{ (\rho,\psi,B):B\ \mbox{is finite dimensional
     C*-algebra,} \\
&&\qquad\qquad \rho:A\to B,\ \psi:B\to A\ \mbox{are unital completely
     positive maps} \}
\end{eqnarray*}
For $\delta>0$, $\omega\in P_f(A)$, and $H\in A_{\sa}$ put
\[
P(H,\omega,\delta)=\inf\{ \log\Tr_B(e^{-\rho(H)}):
(\rho,\psi,B)\in{\rm CPA}(A),\|(\psi\circ\rho)(x)
-x\|<\delta,\forall x\in \omega\}\;,
\]
where the $\inf$ is taken over all $(\rho,\psi,B)\in{\rm CPA}(A)$. Let
\begin{eqnarray*}
&&P_\alpha(H,\omega,\delta)=\limsup_{n\to\infty} \frac{1}{n}
    P\Big( \sum_{j=0}^{n-1}\alpha^j(H),
   \bigcup_{j=0}^{n-1}\alpha^j(\omega),\delta)\;, \\
&&P_\alpha(H,\omega)=\sup_{\delta>0} P_\alpha(H,\omega,\delta)\;.
\end{eqnarray*}

\begin{definition}\label{D11.1}
The {\em pressure of $\alpha$ at\/} $H$ is
\[
P_\alpha(H)=\sup_{\omega\in P_f(A)} P_\alpha(H,\omega)\;.
\]
\end{definition}

Note that $\Tr_B(1)=\rank B$, so when $H=0$ the above definition
restricts to Voiculescu's definition of topological entropy
$ht(\alpha)$. The basic properties of pressure are contained in

\begin{proposition}\label{P11.2}
The following properties are satisfied by $P_\alpha$ for $H,K\in
A_{\sa}$.

(i) \
If $\varphi$ is an $\alpha$-invariant state of $A$ then
\[
P_\alpha(H)\geq h_\varphi(\alpha)-\varphi(H)\;,
\]

\qquad
where $h_\varphi(\alpha)$ is the CNT-entropy of $\alpha$.

(ii) \
If $H\leq K$ then $P_\alpha(H)\geq P_\alpha(K)$.

(iii) \
$P_\alpha(H+c1)=P_\alpha(H)-c$, $c\in\R$.

(iv) \
$P_\alpha(H)$ is either infinite for all $H$ or is finite valued.

(v) \
If $P_\alpha$ is finite valued then $|P_\alpha(H)-P_\alpha(K)|\leq
\|H-K\|$.

(vi) \
For $k\in\N$ $P_{\alpha^k}\big(\sum\limits_{j=0}^{k-1}
\alpha^j(H)\big)=
kP_\alpha(H)$.

(vii) \
$P_\alpha(H+\alpha(K)-K)=P_\alpha(H)$.
\end{proposition}

The proof of (i) is a modification of \cite[Prop.~4.6]{V} using
(\ref{e11.1}). The others are modifications of the
corresponding results in the classical case, \cite[Thm.~9.7]{W}. In the
proof of (ii), (v), (vi) we make use of the important Peierls-Bogoliubov
inequality \cite[Cor.~3.15]{O-P},
\begin{equation}\label{e11.4}
\log\Tr_B(e^H)\leq \log\Tr_B(e^K)\qquad \mbox{if $\;H\leq K$}\;.
\end{equation}
The pressure function has very strong properties, as our next result
shows.

\begin{proposition}\label{P11.3}
Suppose $ht(\alpha)<\infty$, and let $\varphi$ be a self-adjoint linear
functional on $A$. Then $\varphi$ is an $\alpha$-invariant state if and
only if $-\varphi(H)\leq P_\alpha(H)\;$ $\forall H\in A_{\sa}$.
\end{proposition}

The proof is a good illustration of the basic properties of $P_\alpha$.
If $\varphi$ is an $\alpha$-invariant state then by (i) in
Proposition~\ref{P11.2}, $-\varphi(H)\leq
P_\alpha(H)-h_\varphi(\alpha)\leq P_\alpha(H)$.

Conversely, if $-\varphi(H)\leq P_\alpha(H)$ for all $H$ then by (vii)
applied to $H=0$,
\begin{eqnarray*}
-\varphi(\alpha(K)-K) &=& -\frac{1}{n}\varphi
      (\alpha(nK)-nK)\leq\frac{1}{n}P_\alpha(\alpha(nK)-nK) \\
&=& \frac{1}{n}P_\alpha(0)\to0\qquad \mbox{as}\quad n\to\infty\;.
\end{eqnarray*}
Application of this to $-K$ yields $\alpha$-invariance of $\varphi$.

By properties (ii) and (iii)
\begin{eqnarray*}
-\varphi(H) &=& -\frac{1}{n}\varphi(nH)\leq
     \frac{1}{n}P_\alpha(nH) \\
&\leq& \frac{1}{n}P_\alpha(0-n\|H\|)=\frac{1}{n}P_\alpha(0)
+\|H\|\to\|H\|\;.
\end{eqnarray*}
Thus $\|\varphi\|\leq 1$. Since for $c\in\R$,
\[
-c\varphi(1)\leq P_\alpha(c1)=ht(\alpha)-c
\]
we see that $\varphi(1)=1$, so $\varphi$ is a state.
\hfill
$\Box$

\begin{definition}\label{D11.4}
We say an $\alpha$-invariant state is an {\em equilibrium state
at\/} $H$ if
\[
P_\alpha(H)=h_\varphi(\alpha)-\varphi(H)\;,
\]
or equivalently by property (i)
\[
h_\varphi(\alpha)-\varphi(H)=\sup_\psi(h_\psi(\alpha)-\psi(H))\;,
\]
where the $\sup$ is taken over all $\alpha$-invariant states.
\end{definition}

In the finite dimensional case an equilibrium state corresponds to the
Gibbs state.

Recall that if $F$ is a real convex function on a real Banach space
$X$ then a linear functional $f$ on $X$ is called a {\em tangent
functional to $F$ at\/} $x_0\in X$ if
\[
F(x_0+x)-F(x_0)\geq f(x)\qquad (=f(x_0+x)-f(x_0))
\]
for all $x\in X$.

\begin{proposition}\label{P11.5}
Suppose $ht(\alpha)<\infty$ and the pressure is a convex function on
$A_{\sa}$.
\begin{itemize}
\item[(i)]
If $\varphi$ is an equilibrium state at $H$ then $-\varphi$ is a
tangent functional to $P_\alpha$ at $H$.

\item[(ii)]
If $-\varphi$ is a tangent functional for $P_\alpha$ at $H$ then
$\varphi$ is an $\alpha$-invariant state.
\end{itemize}
\end{proposition}

The proof of (i) is immediate from property (i) in
Proposition~\ref{P11.2}. Since
$P_\alpha(H)=h_\varphi(\alpha)-\varphi(H)$,
\[
P_\alpha(H+K)-P_\alpha(H)\geq(h_\varphi(\alpha)-\varphi(H+K))
-(h_\varphi(\alpha)-\varphi(H))=-\varphi(K)\;.
\]
The proof of (ii) is similar to that of Proposition~\ref{P11.3}.
\bigskip

With the definition and the main properties of pressure well taken
care of, we now embark on the variational principle. The setting
will be a restricted class of asymptotically abelian C*-algebras.

\begin{definition}\label{D11.6}
A unital C*-dynamical system $(A,\alpha)$ is called {\em
asymptotically abelian with locality\/} if there is a dense
$\alpha$-invariant $\ast$-subalgebra ${\mathcal A}$ of $A$ such
that for each pair $a,b\in{\mathcal A}$ the C*-algebra they
generate is finite dimensional, and for some $p=p(a,b)\in\N$ we
have $[\alpha^j(a),b]=0$ whenever $|j|\geq p$.
\end{definition}

We call elements of ${\mathcal A}$ {\em local operators\/} and
finite dimensional C*-subalgebras of ${\mathcal A}$ {\em local
algebras}. We may add the identity 1 to ${\mathcal A}$, so we
assume $1\in{\mathcal A}$. Since each finite dimensional C*-algebra
is singly generated an easy induction argument shows
$C^\ast(a_1,\ldots,a_r)$ is a local algebra whenever
$a_1,\ldots,a_r\in{\mathcal A}$. In particular if $A$ is separable
then $A$ is an AF-algebra. A similar argument shows that for each
local algebra $N$ there is $p\in\N$ such that $\alpha^j(N)$
commutes with $N$ whenever $|j|\geq p$.

\begin{theorem}\label{T11.7}
Let $(A,\alpha)$ be a unital separable C*-dynamical system which is
asymptotically abelian with locality. Let $H\in A_{\sa}$. Then
\[
P_\alpha(H)=\sup_\varphi(h_\varphi(\alpha)-\varphi(H))\;,
\]
where the $\sup$ is taken over all $\alpha$-invariant states. In
particular, the topological entropy satisfies
\[
ht(\alpha)=\sup_\varphi h_\varphi(\alpha)\;.
\]
\end{theorem}

For the proof two lemmas are needed. We first consider the simplest
case, which is close to that of shifts on infinite tensor products.
To simplify notation we put for each $k\in\N$, and local algebra $N$
\[
H_k=\sum_{j=0}^k \alpha^j(H)\;,\qquad N_k=C^\ast(\alpha^j(N):0\leq
j\leq k)\;.
\]

\begin{lemma}\label{L11.8}
Let $A,\alpha,H$ be as in Theorem~\ref{T11.7}. Suppose there exists
a local algebra $N$ such that $H\in N$, all $\alpha^j(N)$,
$j\not=0$, commute with $N$, and $C^\ast(\alpha^j(N), j\in\Z)=A$.
Then there exists an $\alpha$-invariant state $\varphi$ such that
\[
P_\alpha(H)=h_\varphi(\alpha)-\varphi(H)=
\lim_{k\to\infty}\frac{1}{k+1}\Tr_{N_k}(e^{-H_k})\;.
\]
\end{lemma}

The lemma is easiest to understand if we note that
\begin{eqnarray*}
&&\frac{1}{k+1}(h_\varphi(\alpha^{k+1})-\varphi(H_k))
    =h_\varphi(\alpha)-\varphi(H) \\
&&\frac{1}{k+1}P_{\alpha^{k+1}}(H_k)=P_\alpha(H)\;,
\end{eqnarray*}
and that the expression inside the limit is $\frac{1}{k+1}$ times the
pressure of $H_k$ in $N_k$ as in equation (\ref{e11.1}).

\begin{lemma}\label{L11.9}
Let $\omega\subset{\mathcal A}$ be a finite set containing $H$. Let $N$
be a local algebra with $\omega\subset N$. Choose $p\in\N$ such that
$\alpha^j(N)$ commutes with $N$ for $|j|\geq p$. Let $k\geq p$, and
$M_k=C^\ast(\alpha^{jk}(N_{k-p}); j\in\Z)$. Then
\[
P_\alpha(H,\omega)\leq\liminf_{k\to\infty}
\frac{1}{k}P_{\alpha^k|M_k}(H_{k-p})\;.
\]
\end{lemma}

Note that $\alpha^{jk}(N_{k-p})$ commutes with $N_{k-p}$ for all
$j\not=0$. Furthermore
\[
M_k=C^\ast(\alpha^{jk+i}(N);1\leq i\leq k-p,j\in\Z)\;,
\]
hence the $\alpha^\ell(N)$'s which appear in $M_k$ are those for which
\[
\ell\not\in\bigcup_{j\in\Z} \{ jk-p+1,jk-p+2,\ldots,jk-1\}\;.
\]

The proof consists of showing that the contribution of these $\ell$'s
is negligible for large $k$.

In order to complete the proof of Theorem~\ref{T11.7} apply
Lemma~\ref{L11.8} to $M_k$ and $\alpha^k$. We get an
$\alpha^k$-invariant state $\psi_k$ on $M_k$ such that
\[
h_{\psi_k}(\alpha^k|_{M_k})-\psi_k(H_{k-p})
=P_{\alpha^k|_{M_k}}(H_{k-p})\;.
\]
By Theorem~\ref{T5.2} and Proposition~\ref{P5.3} we can extend
$\psi_k$ to an $\alpha^k$-invariant state $\widetilde{\varphi}_k$ on
$A$ such that
\[
h_{\widetilde{\varphi}_k}(\alpha^k)\geq h_{\psi_k}(\alpha^k|_{M_k})-1\;.
\]
Put
\[
\varphi_k=\frac{1}{k}\sum_{j=0}^{k-1}
\widetilde{\varphi}_k\circ \alpha^j\;.
\]
Then $\varphi_k$ is $\alpha$-invariant, and by concavity of $\varphi\to
h_\varphi(\alpha)$,
\[
h_{\varphi_k}(\alpha)\geq \frac{1}{k}
h_{\widetilde{\varphi}_k}(\alpha^k)\geq
\frac{1}{k} h_{\psi_k}(\alpha^k|_{M_k})
-\frac{1}{k}\;,
\]
so that
\[
h_{\varphi_k}(\alpha)-\varphi_k(H)\geq
\frac{1}{k}P_{\alpha^k|_{M_k}}(H_{k-p})-\varepsilon\;,
\]
where $\varepsilon$ is small. We then find that
\[
\sup(h_\varphi(\alpha)-\varphi(H))\geq
\liminf_{k\to\infty}\frac{1}{k}
P_{\alpha^k|_{M_k}}(H_{k-p})\;,
\]
and the theorem follows from Lemma~\ref{L11.9}.

\begin{corollary}\label{C11.10}
With our assumptions on $(A,\alpha)$ the pressure is a convex function
of $H$.
\end{corollary}

\begin{proof}
Use the affinity of the function $H\to h_\varphi(\alpha)-\varphi(H)$.

Let $(A,\alpha)$ be asymptotically abelian with locality as before and
$H$ a self-adjoint local operator. Put
\[
\delta_H(x)=\sum_{j\in\Z}[\alpha^j(H),x]\;.
\]
Then $\delta_H$ is a derivation on ${\mathcal A}$ and defines a
one-parameter group $\sigma_t^H=\exp({\rm it}\,\delta_H)$ on $A$. Let
$\beta\geq0$. We say an $\alpha$-invariant state $\varphi$ is an {\em
equilibrium state at  $H$ at inverse temperature\/} $\beta$ if
\[
P_\alpha(\beta H)=h_\varphi(\alpha)-\beta\varphi(H)\quad
(=\sup_\psi(h_\psi(\alpha)-\beta\psi(H)))\;.
\]
\end{proof}

\begin{theorem}\label{T11.11}
Suppose a unital separable C*-dynamical system $({\mathcal A},\alpha)$
is asymptotically abelian with locality, and $ht(\alpha)<\infty$. If
$H$ is a self-adjoint local operator in $A$ and $\varphi$ is an
equilibrium state at $H$ at inverse temperature $\beta$, then
$\varphi$ is a KMS-state for $\sigma^H$ at $\beta$. In particular, if
$ht(\alpha)=h_\varphi(\alpha)$, then $\varphi$ is a trace.
\end{theorem}

In order to prove the theorem we may replace $H$ by $\beta(H)$ and
show $\varphi$ is a KMS-state for $\sigma^H$ at 1. By
Proposition~\ref{P11.5} Theorem~\ref{T11.11} is a consequence of

\begin{theorem}\label{T11.12}
Let $A,\alpha,H$ be as in Theorem~\ref{T11.11}. If $-\varphi$ is a
tangent functional for $P_\alpha$ at $H$ then $\varphi$ is a KMS-state
for $\sigma^H$ at 1.
\end{theorem}

The proof of this theorem is modelled on the corresponding proof for
spin lattice systems, see \cite[Ch.~6]{B-R}.

Several examples encountered in the previous sections are
asymptotically abelian with locality. They are:
\begin{eqnarray}
&&\mbox{Shifts on $\;\bigotimes\limits_{i\in\Z}B_i$, where $B_i=B_0$ is
    an AF-algebra.}
   \label{e11.13} \\ [1ex]
&&\mbox{Binary shifts $\;(A(X),\alpha)$ as defined in
    Section~9 with $X\subset\N$ a finite set,} \label{e11.14} \\
&&\mbox{ see \cite{G-S1}.}
   \nonumber \\
&&\mbox{The shift on the Jones projections in the Temperley-Lieb
    algebra, see~10.4.}
   \label{e11.15} \\
&&\mbox{The canonical shift defined by a subfactor $N\subset M$ of
   finite index, see }
   \label{e11.16} \\
&&\mbox{Section~10.} \nonumber \\
&&\mbox{Let $A=K(H)+\C1$, where $K(H)$ is the compact operators on a
     separable}
    \label{e11.17} \\
&&\mbox{Hilbert space. Let $\alpha=\Ad u|_A$, where
$U$ is the
bilateral shift on $H$}.
\nonumber
\end{eqnarray}

\noindent
{\bf 11.13. Counter examples}\
Theorem~\ref{T11.7} is false without the assumption of asymptotic
abelianness. Indeed Theorem~\ref{T9.9} provides a counter example. In
that case $\tau$ is the unique $\alpha$-invariant state, and
$ht(\alpha)\geq\frac{1}{2}\log 2$.

The assumption $ht(\alpha)<\infty$ is necessary in
Theorem~\ref{T11.11}, as is immediate from example
(11.3)
above. If
$B_0$ is infinite dimensional there exist many $\alpha$-invariant
states with infinite entropy.

Locality is necessary in Theorem~\ref{T11.11}. For let $U$ be a unitary
operator on a separable Hilbert space $H$ with singular spectrum such
that $(U^nf,g)\to0$ as $n\to\infty$, for all $f,g\in H$. Let $A$ be the
even CAR-algebra, i.e. the fixed points in the CAR-algebra of the
Bogoliubov automorphism $\alpha_{-1}$, see Section~4. From the proof of
Lemma~\ref{L4.1} we have $ht(\alpha_U)=0$. However, there exist many
$\alpha_U$-invariant states, e.g. all quasi-free states
$\omega_\lambda$, $0\leq\lambda\leq 1$.


\begin{thebibliography}{9999}
\bibitem[A-F]{A-F}
{\bf Alicki, R.} and {\bf Fannes, M.}
Defining quantum dynamical entropy,
Lett. Math. Phys. 32 (1994), 75--82.

\bibitem[A]{A}
{\bf Araki, H.}
Relative entropy for states of von Neumann algebras II,
Publ. RIMS Kyoto Univ. 13 (1977), 173--192.

\bibitem[Ar]{Ar}
{\bf Arveson, W.}
Subalgebras of C*-algebras,
Acta Math. 123 (1969), 141--224.

\bibitem[Av]{Av}
{\bf Avitzour, D.}
Free products of C*-algebras and von Neumann algebras,
Trans. Amer. Math. Soc. 217 (1982), 423--435.

\bibitem[B]{B}
{\bf Besson, O.}
On the entropy of quantum Markov states,
Lecture Notes in Math. 1136 (1985), 81--89, Springer-Verlag.

\bibitem[B-G]{B-G}
{\bf Bezuglyi, S.I.} and {\bf Golodets, V.Ya.}
Dynamical entropy for Bogoliubov actions of free abelian groups on the
CAR-algebra,
Ergod. Th. \& Dynam. Sys. 17 (1997), 757--782.

\bibitem[Bo-Go]{Bo-Go}
{\bf Boca, F.} and {\bf Goldstein, P.}
Topological entropy for the canonical endomorphism of Cuntz-Krieger
algebras, preprint 1999.

\bibitem[BKRS]{BKRS}
{\bf Bratteli, O., Kishimoto, A., R\o{}rdam, M.} and {\bf St\o{}rmer,
E.}  The crossed product of a UHF-algebra by a shift,
Ergod. Th. \& Dynam. Sys. 13 (1993), 615--626.

\bibitem[B-R]{B-R}
{\bf Bratteli, O.} and {\bf Robinson, D.W.}
Operator algebras and quantum statistical mechanics II.
Springer-Verlag 1981.


\bibitem[Br1]{Br1}
{\bf Brown, N.}
Topological entropy in exact C*-algebras,
Math. Ann. 314 (1999), 347--367.

\bibitem[Br2]{Br2}
{\bf Brown, N.}
A note on topological entropy, embeddings and unitaries in nuclear
quasidiagonal C*-algebras, preprint 1998.

\bibitem[B-C]{B-C}
{\bf Brown, N.} and {\bf Choda, M.}
Approximation entropies in crossed products with an application to free
shifts, preprint 1999.

\bibitem[C1]{C1}
{\bf Choda, M.}
Entropy for *-endomorphisms and relative entropy for subalgebras,
J. Operator Th. 25 (1991), 125--140.

\bibitem[C2]{C2}
{\bf Choda, M.}
Entropy for canonical shifts, Trans. Amer. Math. Soc. 334 (1992),
827--849.

\bibitem[C3]{C3}
{\bf Choda, M.}
Reduced free products of completely postitve maps and entropy for
free products of automorphismes, Publ. RIMS Kyoto Univ. 32 (1996),
371--382.

\bibitem[C4]{C4}
{\bf Choda, M.}
Entropy of Cuntz's canonical endomorhism, Pacific J. Math. (1999),
235--245.

\bibitem[C5]{C5}
{\bf Choda, M.}
A dynamical entropy and applications to canonical endomorphisms,
preprint 1998.

\bibitem[C6]{C6}
{\bf Choda, M.}
Entropy on crossed products and entropy on free products, preprint
1999.

\bibitem[C-H]{C-H}
{\bf Choda, M.} and {\bf Hiai, F.}
Entropy for canonical shifts II. Publ. RIMS Kyoto Univ. 27 (1991),
461--489.

\bibitem[C-N]{C-N}
{\bf Choda, M.} and {\bf Natsume, T.}
Reduced C*-crossed products by free shifts, Ergod. Th. \& Dynam. Sys. 18
(1998), 1075--1096.

\bibitem[Co]{Co}
{\bf Connes, A.}
Entropie de Kolmogorff Sinai et mechanique statistique quantique,
C. R. Acad. Sci. Paris 301 (1985), 1--6.

\bibitem[CFW]{CFW}
{\bf Connes, A.}, {\bf Feldmann, J.} and {\bf Weiss, B.}
An amenable equivalence relation is generated by a single
transformation, Ergod. Th. \& Dynam. Sys. 1 (1981), 431--450.

\bibitem[CNT]{CNT}
{\bf Connes, A.}, {\bf Narnhofer, H.} and {\bf Thirring, W.}
Dynamical entropy of C*-algebras and von Neumann algebras, Commun.
Math. Phys. 112 (1987), 691--719.

\bibitem[C-S]{C-S}
{\bf Connes, A.} and {\bf St\o{}rmer, E.}
Entropy of automorphisms of II$_1$-von Neumann algebras, Acta Math. 134
(1975), 289--306.

\bibitem[Cu]{Cu}
{\bf Cuntz, J.}
Simple C*-algebras generated by isometries, Commun. Math. Phys. 57
(1977), 173--185.

\bibitem[De]{De}
{\bf Deaconu, V.}
Entropy estimates for some C*-endomorphisms, Proc. Amer. Math. Soc. 127
(1999), 3653--3658.

\bibitem[D1]{D1}
{\bf Dykema, K.}
Exactness of reduced amalgamated free products of C*-algebras, preprint
1999.

\bibitem[D2]{D2}
{\bf Dykema, K.}
Topological entropy of some automorphisms of reduced amalgamated free
product C*-algebras, preprint 1999.

\bibitem[D-S]{D-S}
{\bf Dykema, K.} and {\bf Shlyakhtenko, D.}
Exactness of Cuntz-Pimsner C*-algebras, preprint 1999.

\bibitem[E]{E}
{\bf Enomoto, M.}, {\bf Nagisa M.}, {\bf Watatani, Y.} and
{\bf Yoshida, H.}
Relative commutant algebras of Powers' binary shifts on the hyperfinite
II$_1$-factor, Math. Scand. 68 (1991), 115--130.

\bibitem[G-N1]{G-N1}
{\bf Golodets, V.~Ya.} and {\bf Neshveyev, S.}
Dynamical entropy for Bogoliubov actions of torsion-free abelian groups
on the CAR-algebra, preprint 1998.

\bibitem[G-N2]{G-N2}
{\bf Golodets, V.~Ya.} and {\bf Neshveyev, S.}
Entropy of automorphisnes of II$_1$-factors arising from the dynamical
systems theory, preprint 1999.

\bibitem[G-N3]{G-N3}
{\bf Golodets, V. Ya.} and {\bf Neshveyev, S.}
Non-Bernoullian quantum $K$-systems.
Commun. Math. Phys. 195 (1998), 213--232.

\bibitem[G-S1]{G-S1}
{\bf Golodets, V. Ya.} and {\bf St\o{}rmer, E.}
Entropy of C*-dynamical systems defined by bitstreams, Ergod. Th. \&
Dynam. Sys. 18 (1998), 859--874.

\bibitem[G-S2]{G-S2}
{\bf Golodets, V. Ya.} and {\bf St\o{}rmer, E.}
Generators and comparison of entropies of automorphisms of finite von
Neumann algebras, J. Funct. Anal. 164 (1999), 110--133.

\bibitem[H-S]{H-S}
{\bf Haagerup, U.} and {\bf St\o{}rmer, E.}
Maximality of entropy in finite von Neumann algebras, Invent. Math. 132
(1998), 433--455.

\bibitem[H]{H}
{\bf Hiai, F.}
Entropy for canonical shifts and strong amenability,
Int. J. Math.
6 (1995), 381--396.

\bibitem[Hu]{Hu}
{\bf Hudetz, T.}
Quantum dynamical entropy revisited,
Quantum probability (1997), 241--251, Banach Center Publ. 43.

\bibitem[J]{J}
{\bf Jones, V.F.R.}
Index for subfactors, Invent. Math. 72 (1983), 1--25.

\bibitem[K]{K}
{\bf Kirchberg, E.}
On subalgebras of the CAR-algebra, J. Funct. Anal. 129 (1995), 35--63.

\bibitem[L]{L}
{\bf Longo, R.}
Index of subfactors and statistics of quantum fields, Commun. Math.
Phys. 126 (1989), 145--155.

\bibitem[M]{M}
{\bf Moriya, H.}
Variational principle and the dynamical entropy of space translations,
preprint 1999.

\bibitem[N-U]{N-U}
{\bf Nakamura, M.} and {\bf Umegaki, H.}
A note on the entropy for operator algebras, Proc. Japan Acad. 37
(1961), 149--154.

\bibitem[NST]{NST}
{\bf Narnhofer, H., St\o{}rmer, E.} and {\bf Thirring, W.}
C*-dynamical systems for which the tensor product formula for entropy
fails,
Ergod. Th. \& Dynam. Sys. 15 (1995), 96/9/68.

\bibitem[N-T1]{N-T1}
{\bf Narnhofer, H.} and {\bf Thirring, W.}
Dynamical entropy of quantum systems and their abelian counterpart. On
Klauder Path: A field trip, Emch, G.G., Hegerfeldt, G.C., Streit, L.
World Scientific; Singapore, 1994.

\bibitem[N-T2]{N-T2}
{\bf Narnhofer, H.} and {\bf Thirring, W.}
C*-dynamical systems that are highly anticommutative,

\bibitem[N]{N}
{\bf Neshveyev, S.}
Entropy of Bogoliubov automorphisnes of CAR and CCR algebras with
respect to quasi-free states,Rev. Math. Phys.

\bibitem[N-S1]{N-S1}
{\bf Neshveyev, S.} and {\bf St\o{}rmer, E.}
Entropy in type I algebras, preprint 2000.

\bibitem[N-S2]{N-S2}
{\bf Neshveyev, S.} and {\bf St\o{}rmer, E.}
The variational principle for a class of asymptotically abelian
C*-algebras, preprint 2000.

\bibitem[O-P]{O-P}
{\bf Ohya, M.} and {\bf Petz, D.}
Quantum entropy and its use. Texts and Monographs in Physics,
Springer-Verlag, 1993.

\bibitem[O]{O}
{\bf Ornstein, D.S.}
Bernoulli shifts with the same entropy are isomorphic, Adv. Math. 4
(1970), 337--352.

\bibitem[P-S]{P-S}
{\bf Park, Y.M.} and {\bf Shin, H.H.}
Dynamical entropy of space translations of CAR and CCR algebras with
respect to quasi-free states,
Commun. Math. Phys. 152 (1993), 497--537.

\bibitem[P-P]{P-P}
{\bf Pimsner, M.} and {\bf Popa, S.}
Entropy and index for subfactors, Ann. Sci. Ecole Norm. Sup. 19 (1986),
57--106.

\bibitem[PWY]{PWY}
{\bf Pinzari, C., Watatani, Y.} and {\bf Yonetani, K.}
KMS states, entropy and the variational principle in full C*-dynamical
systems, preprint.

\bibitem[P]{P}
{\bf Powers, R.T.}
An index theory for semigroups of *-endomorphisms of B(H) and type
II$_1$-factors, Canad. J. Math. 40 (1988), 86--114.

\bibitem[Po-Pr]{Po-Pr}
{\bf Powers, R.T.} and {\bf Price, G.}
Binary shifts on the hyperfinite II$_1$-factor, Contemp. Math. 145
(1993), 453--464.

\bibitem[Pr]{Pr}
{\bf Price, G.}
The entopy of rational Powers shifts, Proc. Amer. Math. Soc.

\bibitem[Q]{Q}
{\bf Quasthoff, V.}
Shift automorphisms of the hyperfinite factor, Math. Nachrichten 131
(1987), 101--106.

\bibitem[Sa]{Sa}
{\bf Sauvageot, J.}
Ergodic properties of the action of a matrix in $SL(2,\Z)$ on a non
communative torus, preprint.

\bibitem[S-T]{S-T}
{\bf Sauvageot, J.-L.} and {\bf Thouvenot, J.-P.}
Une nouvelle definition de l'entropie dynamique des systems non
commutatifs, Commun. Math. Phys. 145 (1998), 521--542.

\bibitem[Sh]{Sh}
{\bf Shields, P.}
The theory of Bernoulli shifts, Univ. Chicago Press, 1973.

\bibitem[S-S]{S-S}
{\bf Sinclair, A.} and {\bf Smith, R.R.}
The completely bounded approximation property for discrete crossed
products, Indiana Univ. Math. J. 46 (1997), 1311--1322.

\bibitem[S1]{S1}
{\bf St\o{}rmer, E.}
Entropy of some automorphisms of the II$_1$-factor of the free group in
infinite number of generators, Invent. Math. 110 (1992), 63--73.

\bibitem[S2]{S2}
{\bf St\o{}rmer, E.}
Entropy of some inner automorphisms of the hyperfinite II$_1$-factor,
Int. J. Math. 4 (1993), 319--322.

\bibitem[S3]{S3}
{\bf St\o{}rmer, E.}
States and shifts on infinite free products of C*-algebras, Fields Inst.
Commun. 12 (1997), 281--291.

\bibitem[S4]{S4}
{\bf St\o{}rmer, E.}
Entropy of endomorphisms and relative entropy in finite von Neumann
algebras, J. Funct, Anal. 171 (2000), 34--52.

\bibitem[SV]{SV}
{\bf St\o{}rmer, E.} and {\bf Voiculescu, D.}
Entropy of Bogoliubov automorphisms of the canonical anticommutation
relations, Commun. Math. Phys. 133 (1990), 521--542.

\bibitem[T]{T}
{\bf Thomsen, K.}
Topological entropy for endomorphisms of local C*-algebras,
Commun. Math. Phys. 164 (1994), 181--193.

\bibitem[Vi]{Vi}
{\bf Vik, S.}
Fock representation of the binary shift algebra, Math. Scand.

\bibitem[V]{V}
{\bf Voiculescu, D.}
Dynamical approximation entropies and topological entropy in operator
algebras, Commun. Math. Phys. 170 (1995), 249--281.

\bibitem[W]{W}
{\bf Walters, P.}
An introduction to ergodic theory, Graduate texts in Math. 79,
Springer-Verlag, 1982.
\end{thebibliography}
\end{document}